\documentclass [11pt]{article}
\usepackage{amsmath, amsthm, amssymb}
\usepackage{graphicx}
\newfont{\bbb} {msbm10}
\newcommand{\R}{\Bbb{R}}
\newcommand{\bS}{\Bbb{S}}

\newcommand{\T}{\Bbb{T}}

\newcommand{\Q}{\Bbb{Q}}
\newcommand{\D}{\Bbb{D}}

\newcommand{\sbs}{\subset}
\newcommand{\ra}{\rightarrow}

\newcommand{\dsigma}{{\dot{\sigma}}}

\newcommand{\p}{\partial}

\newcommand{\cL}{{\cal{L}}}

\newcommand{\cS}{{\cal{S}}}

\newcommand{\cA}{{\cal{A}}}

\newcommand{\HH}{\Bbb{H}}

\newcommand{\s}[1]{{\sf{#1}}}

\newcommand{\0}[1]{_{_{#1}}}
\newcommand{\rC}{{\rm{C}}\,}

\newcommand{\sN}{{\sf{N}}}
\newcommand{\sL}{{\sf{Link}}}

\newcommand{\dX}{\dot{X}}

\newcommand{\dsquare}{\dot{\square}}
\usepackage[all]{xy}

\begin{document}

\title{Normal Smoothings for Charney-Davis Strict Hyperbolizations}
\author{Pedro Ontaneda\thanks{The author was
partially supported by a NSF grant.}}
\date{}

\maketitle

\begin{abstract} 
We prove that the Charney-Davis Strict
Hyperbolization of a smoothly cubulated manifold
admits a {\it normal smooth structure}.
We also prove that this normal smooth structure
is diffeomorphic to a 
smooth structure that has good tangential properties.

The results in this paper are key ingredients in problem of smoothing the metric of a strictly hyperbolized manifold (see  \cite{ChD}, \cite{O2}). 

\end{abstract}\vspace{.2in}

In his 1987 paper ``Hyperbolic Groups''
\cite{G} M. Gromov introduced 
the process of {\it hyperbolization}.
This process assigns to each simplicial complex $K$
a nonpositively curved (in the geodesic sense) complex.
The hyperbolization process has a lego type flavor and 
it can roughly be described in simple terms: to construct 
the hyperbolization of a simplicial complex $K$
we replace its basic set of pieces (simplices)
by another
basic set of ``hyperbolization pieces". In other
words,  to construct the hyperbolization of $K$
we assemble the hyperbolization pieces
using the same pattern as the one used to assemble
$K$. 
The hyperbolization process was later
studied and used by Davis and  Januszkiewicz in \cite{DJ}.\\

An important property of hyperbolization is
that if $K$ is a $PL$ manifold, then the hyperbolization of 
$K$ is also a $PL$ manifold. Moreover, Davis and 
Januszkiewicz \cite{DJ}
showed that the hyperbolization of a smoothly triangulated
manifold is also a smooth manifold.\\

%if the complex $K$ is a smooth
%manifold (the smooth structure is compatible
%with $K$) then the hyperbolization of $K$ is
%also a smooth manifold.\\ 

In \cite{ChD}  Charney and Davis 
built on previous versions of hyperbolization
and presented the {\it strict hyperbolization} process.
In this case one begins with a cube complex $K$
(with large links) and obtains a negatively curved
space $K_X$. In this process we replace the cubes by
what we call {\it Charney-Davis strict hyperbolization pieces}.
Again,  the hyperbolization of a smoothly cubulated manifold is also a smooth manifold.\\
%That is, if $K$ is smoothable, then so is $K_X$.\\

Notation: A smooth cubulation of a smooth manifold
$M$ is a homeomorphism $f:K\ra M$, where
$K$ is a cube complex and $f$ is a smooth embedding
when restricted to each cube of $K$.
In this case, for simplicity, we will just say that {\it $K$ is a smooth cube complex}, or a {\it smooth cube manifold}. Therefore, if $K$ is a smooth cube complex, then $K_X$ is smoothable.
\\

%if $K$ is a smooth
%manifold then the strict hyperbolization of $K$ is
%also a smooth manifold.\\

Let $K$ be a smooth cube complex. Then the
Charney-Davis-Januszkiewicz smooth structure on $K_X$
is ``good" from the Geometric Topology point of view
because it has good tangential properties
(see Section 4). But it is
quite poor from the Geometry point of view
because there is no  a priori relationship between
the smooth structure and the rich cube geometry of $K_X$.
To correct this
we introduce ``normal smooth structures"
in the paragraphs below. These normal structures are very natural
and useful (see \cite{O2}). Of course, there is no guarantee
that such structures exist on $K_X$; actually there is
no guarantee that $K_X$ is smoothable at all (see \cite{BLW}
for examples of non-smoothable $K_X$).
In this paper we prove the following.\\

\noindent {\bf Main Theorem.}
{\it Let $K$ be a cube complex. If $K$ is smooth, then $K_X$ admits a normal smooth structure.}\\

\noindent {\bf Addendum.} {\it  Furthermore, the normal smooth structure
has good tangential properties}.\\
 
By  ``good tangential properties" we mean that
$K_X$, with its normal smooth structure, smoothly embeds
in $K\times X$ with trivial normal bundle. (Here the cube complex
$K$ is considered with a ``normal smooth structure" as in
\cite{O1}, see also Theorem 4.3
for more details). In \cite{ChD} Charney and Davis ask the
following:\\

\noindent{\it Question 7.4 in \cite{ChD}. Is there a 
stably parallelizable Charney-Davis strict hyperbolization
piece?}\\

The addendum together with a positive answer to
the question above would imply that if $K$ is a smooth
cube manifold then the natural map
$K_X\ra K$ is covered by a map of stable tangent bundles,
where we are considering $K_X$ here with a normal smooth structure.\\

 Before we explain what normal structures are,
 we recall some features of the
 Charney-Davis strict hyperbolization process.
 For more details see Section 3.\\

 We write $\square^n=[0,1]^n$.
A {\it Charney-Davis hyperbolization piece $X^n$ of dimension $n$} is a compact
connected orientable hyperbolic $n$-manifold with corners satisfying the properties stated in Lemma 6.2 of \cite{ChD}.
We state some relevant facts.
There is a smooth map $f:X^n\ra\square^n$, such that
$X^n$ and $f$ satisfy the following.

\begin{enumerate}
\item[(1)] For any $k$-face $\square^k$ of $\square^n$ we have that the $k$-face $X_{\square^k}=f^{-1}
(\square^k)$ of $X^n$ is totally geodesic in $X^n$, and it is a Charney-Davis hyperbolization
piece of dimension $k$. 

\item[(2)] The faces of $X^n$ intersect orthogonally
(unless one is a face of the other).
\end{enumerate}

\noindent The interior $X_{\square^k}=f^{-1}(\dot{\square^k})$ will be denoted by $\dX_{\square^k}$.\\

Let $X_{\square^k}$ be a $k$-face. We denote by 
$\sL(X_{\square^k},X^n)$ the link 
of $X_{\square^k}$ in $X^n$ (at $p$), that is, the set of inward vectors orthogonal to
$X_{\square^k}$ at $p$, for some $p\in X_{\square^k}$. The
link can be identified with the canonical all-right spherical $(n-k-1)$-simplex
$\Delta^{n-k-1}$. In this sense we consider $\Delta^{n-k-1}\sbs T_pX$.
Similarly we can consider the link $\sL(\square^k,\square^n)$
of $\square^k$ in $\square^n$. It can also be
identified with $\Delta^{n-k-1}\sbs T_q\square^n$, for 
$q\in \square^k$. 
We can identify the derivative of $f$,
$(Df_p|\0{\Delta^{n-k-1}}):\Delta^{n-k-1}
\ra\Delta^{n-k-1}$,
with the identity map.\\

Let $K$ be a cube complex and $\square\in K$.
Recall that the link $\sL(\square,K)$
of $\square$ in $K$ is the all-right spherical complex
$\{\sL(\square,\square')\,:\, \square\sbs\square'\in K\}$.\\

The strict hyperbolization process of Charney and Davis is done by gluing copies
of $X^n$ using the same pattern as the one used to obtain the cube complex $K$
from its cubes (see Section 3 for more details). This space is called $K_X$ in \cite{ChD}. 
Note that we get a map
$F:K_X\ra K$, which restricted to each copy of $X^n$ is just the map $f:X^n\ra\square^n$. We will write $X_{\square^k}=F^{-1}(\square^k)$, for a $k$-cube $\square^k$
of $K$.\\

The link $\sL(X_{\square},X_K)$
of $X_{\square}$ in $X_K$ is the all-right spherical complex
$\{\sL(X_{\square},X_{\square'})\,:\, \square\sbs\square'
\in K\}$.\\

We can use the derivative of
the map $F:K_X\ra K$ (in a piecewise fashion) to identify $\sL(X_{\square^k},K_X)$ with $\sL(\square^k,K)$. Hence we write
$\sL(X_{\square^k},K_X)=\sL(\square^k,K)$; thus the set of links of faces $\square$ in $K$ coincides with the set of links 
of the $X_\square$ in $K_X$.\\

Let $\square\in K$. A {\it link smoothing} of 
$\sL(\square^k,K)=\sL(X_{\square^k},K_X)$
is just a homeomorphism
$h_{\square^k}:\bS^{n-k-1}\ra\sL(X_{\square^k},K_X)$. A {\it (complete) set of link
smoothings} is a set $\cA=\{h_\square\}\0{\square\in K}$. 
\\

We are now ready to define normal structures.
Let $X_{\square^k}\sbs X_{\square^{n}}$ be a $k$-face of $K_X$,
contained in the copy $X_{\square^n}$ of $X$ over $\square^n$. For a non-zero vector $u$ normal to $X_{\square^k}$ at $p\in X_{\square^k}$, and pointing inside 
$X_{\square^n}$, we have that $exp_p(tu)$ is defined and contained in
$X_{\square^n}$, for $0\leq t< t_0$,
for certain fixed $t_0$. Let $\cA=\{h_\square\}\0{\square\in K}$ be a fixed set of link smoothings.
We define the map 
$$H\0{\square^k}\,\,\,:\,\,\D^{n-k}\times \dX_{\square^k}
\,\,\,\,\longrightarrow\,\,\,\, K_X $$
\noindent given by
$$H\0{\square^k}(\,t\,v\,\, ,\,\,p\,)\,\,\,=\,\,\, exp\0{p}\,\Big(\,\,2r\,t\,\,h_{\square^k}(v)\,\,\Big)$$
\noindent where $v\in\bS^{n-k-1}$ and $t\in [0,1)$.
We have that $H\0{\square^k}$ is a topological embedding. 
The map $H\0{\square^k}$ is called a {\it normal chart for the $k$-face $X_{\square^k}$}.
The collection $\big\{ H\0{\square^k} \big\}\0{\square^k\in K}$ of the normal charts
is a {\it normal atlas}, and if this atlas is smooth (or $C^k$) the induced
differentiable structure is called a {\it normal smooth (or $C^k$) structure.}\\

Note that the charts $H\0{\square}$ respect normal
directions and radial distances. Note also that
the normal atlas  
$\big\{ H\0{\square^k} \big\}\0{\square^k\in K}$
depends only on the set $\cA$ of links smoothings.\\

Here is a description of the paper. In Section 1 we deal with
smooth structures on cube complexes. In this section we
recall and review some necessary concepts, definitions and results that appear in \cite{O1}. In Section 2 we
study Charney-Davis strict hyperbolization pieces.
In Section 3 we review the Charney-Davis strict
hyperbolization process, compare two ways of doing this
process. In this section we also introduce ways to construct
smooth structures on Charney-Davis hyperbolizations
with good tangential properties. In Section 4 we
deal with normal smooth structures. We also treat
the case of smooth manifolds with one point
singularities. There are three appendices.\\

The results in this paper are key ingredients in problem of smoothing the metric of a strictly hyperbolized manifold (see \cite{O2}).

\vspace{.4in}

\begin{center} {\bf \large 1. Smooth Structures on Cube Complexes
and All-Right
Spherical Complexes.}
\end{center}

\noindent {\bf 1.1. Normal Smooth Structures.}\\

For the basic definitions and results about cube and spherical  complexes see
for instance \cite{BH}. Recall that a spherical
complex is an {\it all-right spherical complex} if all of its edge lengths are equal to 
$\pi/2$. 
Given a (cube or all-right spherical) complex $K$ we use the same notation $K$ for the complex itself (the collection of all
closed cubes or simplices) and its realization (the union of all cubes or simplices).
For $\sigma\in K$ we denote its interior by $\dsigma$.\\

Let $M^n$ be a smooth manifold of dimension $n$. A {\it smooth cubulation} of  $M$ is given by $(K,f)$, where 
$K$  a cube complex and $f:K\ra M$ a non-degenerate
piecewise differentiable homeomorphism \cite{MunkresLectures}, that is, for all $\sigma\in K$  we have $f|_{\sigma}$ is a smooth embedding. 
Sometimes we will write $K$ instead
of $(K,f)$. The smooth manifold $M$ together with a smooth cubulation is
a {\it smooth cube manifold} or a {\it smooth cube complex}. A {\it smooth all-right-spherical triangulation}  and a 
{\it smooth all-right-spherical manifold} (or complex) is defined analogously.\\

Note that if $K$ is a smooth cubulation (or all-right
spherical triangulation) of $M$, then $K\cong_{PL}M$, that is,
$K$ is $PL$-homeomorphic to the smooth manifold $M$.\\

The {\it geometric link} 
 $\sL(\sigma^{j},K)$ of an open $j$-cube or $j$-all-right simplex $\sigma^j$
is the union of the end points of straight (geodesic) segments of small length $\epsilon>0$
emanating perpendicularly (to $\dsigma^j$) from some point  $x\in \dsigma^j$. 
We say that the link is {\it based at $x$.} 
And the star $\s{Star}(\sigma,K)$ as the union of such segments. We can identify the star with the cone of the link
$\rC\sL (\sigma,K)$  (or $\epsilon$-cone) defined as
$$\rC \, \sL(\sigma, K)=\sL(\sigma,K)\times [0,\epsilon)\,/\, \sL(\sigma,K)\times\{0\}.$$ 
\noindent  We shall denote the {\it cone point} by $o$ or, more specifically, by 
$o\0{\rC \, \sL(\sigma, K)}$.
Thus a point $x$ in $\rC \, \sL(\sigma, K)$, different from the cone point $o$, can be written as
$x=t\,u$, $t\in (0,\epsilon)$, $u\in \sL(\sigma, K)$. For $s>0$ we get the
{\it cone homothety} $x\mapsto sx=(st)u$ (partially defined if $s>1$). 
If we want to make explicit the dependence of the link or the cone on $\epsilon$ we shall write $\sL_\epsilon(\sigma,K)$ or $\rC _\epsilon\,\sL(\sigma,K)$ respectively. 
Also, we will always take $\epsilon <1/2$ ($<\pi/4$ in the spherical case) and it can be verified
that all results in this section (unless otherwise stated) are independent
of the choice of the $\epsilon$'s. 
As usual we shall identify the $\epsilon$-neighborhood
of $\dsigma$ in $K$ with $\rC _\epsilon\,\sL(\sigma,K)\times
\dsigma$ (or just $\rC \,\sL(\sigma,K)\times
\dsigma$). Hence a cone homothety induces a {\it neighborhood homothety}
obtained by crossing it with the identity $1\0{\dsigma}$.
Note  that $\sL(\sigma,K)$ and $\rC \,\sL(\sigma,K)$ are subsets of $K$. 
\\

In what follows we assume that $f:K\ra M$ is a smooth cubulation
(or all-right spherical triangulation) of the smooth manifold $M$.
 Recall that the link $\sL(\sigma^i,K)$,  $\sigma^i\in K$, has a natural
all-right piecewise spherical structure, which induces a simplicial structure and thus a $PL$ structure on  $\sL(\sigma^i,K)$.
Since the $PL$ structure on $M$ induced by $K$ is Whitehead compatible with $M$ we have that
the link $\sL(\sigma^i,K)$ is $PL$ homeomorphic to $\bS^{n-i-1}$.
A {\it link smoothing for}  $\dsigma^i$ (or $\sigma^i$)
is just a homeomorphism $h_{\sigma^i}:\bS^{n-i-1}\ra\sL(\sigma^i,K)$.
The {\it cone} of $h_{\sigma^i}$ is the map 
$$\rC \,h_{\sigma^i}:\D^{n-i}\longrightarrow \rC\sL (\sigma^i, K)$$
\noindent given by $t\,x= [x,t]\mapsto t\,h_{q^i}(x)=[h_{q^i}(x),  \,t]$,
 where we are canonically identifying the $\epsilon$-cone  of $\bS^{n-i-1}$ with the disc $\D^{n-i}$. 
 We remark that we are not assuming $h_{\sigma^i}$ to be smooth (or even $PL$).\\

A link smoothing $h_{\sigma^i}$  induces the following  smoothing of the normal neighborhood of
$\dsigma^i$:
$$h^\bullet_{\sigma^i}=f\,\, \circ\,\,\Big(\rC \,h_{\sigma^i}\times 1_{\dsigma^i}\Big):\D^{n-i}\times \dsigma^i\longrightarrow M. $$

The pair $(\,h^\bullet_{\sigma^i}\, ,\,\D^{n-i}\times \dsigma^i\,)$, or simply 
$h^\bullet_{\sigma^i}$, is a {\it normal
chart} on $M$. Note that the collection  $\cA=\big\{\,(\,h^\bullet_{\sigma^i}\, ,\,\D^{n-i}\times 
\dsigma^i\,)\,\big\}_{\sigma^i\in K}$ is a topological atlas for $M$.
Sometimes will just write $\cA=\big\{\,h^\bullet_{\sigma^i}\,\big\}_{\sigma^i\in K}$.
The topological atlas $\cA$ is called a {\it normal atlas}. It depends uniquely on 
the complex $K$, the map $f$ and the collection of link smoothings 
$\{h_{\sigma}\}_{\sigma\in K}$. 
To express the dependence of the atlas on the set of links smoothings
we shall write $\cA=\cA\big(\{h_{\sigma}\}_{\sigma\in K}\,\big)$
(this is different from $\cA=\big\{\,h^\bullet_{\sigma^i}\,\big\}_{\sigma^i\in K}$,
as written above).\\

The most important feature about these normal
atlases is that they preserve the radial and sphere (link) structure given by $K$.
These features make normal atlases very powerful tools for geometric
constructions.\\

Note that not every collection of link smoothings induce a smooth atlas. But when 
the atlas is smooth, we call $\cA$ a {\it normal smooth atlas on $M$ with respect to} $K$ and the corresponding smooth structure $\cS'$
a {\it normal smooth structure on $M$ with respect to $K$}. \\

\noindent {\bf Remarks.}

\noindent {\bf 1.} If the normal atlas $\cA$  is smooth the maps $f|_{\dsigma^i}:\dsigma^i\ra (M,\cS')$ and the link smoothings
$$h_{\sigma^i}:\bS^{n-i-1}\ra (M,\cS')$$ 
\noindent are, by construction, smooth embeddings. Here, as before,
$\cS'$ is the normal smooth structure induced by the normal smooth atlas  $\cA$.

\noindent {\bf 2.} The atlas $\cA$ is smooth if and only if there is a smooth structure
$\cS'$ such that all normal charts 
$h^\bullet_{\sigma^i}:\D^{n-i}\times \dsigma^i\longrightarrow (M,\cS') $
are smooth embeddings. (This is true for any topological atlas.)\\

Later in this section we  show that the
atlas $\cA\big(\{h_{\sigma_i}\}\big)$ is smooth if and only if
the set of link smoothings $\{h_{\sigma_i}\}$ is ``smoothly compatible". Here is the main result of \cite{O1}.\\

\noindent {\bf Theorem 1.1.1.} {\it Let $M$ be a smooth cube manifold,
with smooth structure $\cS$.
Then $M$ admits a normal smooth structure $\cS'$ diffeomorphic to $\cS$.}\\

Hence if $M^n$ is a smooth manifold with smooth structure $\cS$ and $K$ is a cubulation of $M$, then there are link smoothings $ h_{\sigma^i}$, for all $\sigma^i\in K$, such that the atlas $\cA=\cA\big(\{h_{\sigma}\}_{\sigma\in K}\,\big)$
is smooth. Moreover the normal smooth structure $\cS'$,
induced by $\cA$, is diffeomorphic to $\cS$.\\

\noindent {\bf Addendum Theorem 1.1.1.} {\it The statement of Theorem 1.1.1 also holds for smooth all-right-spherical complexes.}\\ 

The following is a corollary of the proof of Theorem 1.1.1 given
in \cite{O1} (see Lemma 1.2 in \cite{O1}).\\

\noindent {\bf Corollary 1.1.2.} {\it Let $f:K\ra (M,\cS)$
be a smooth cubulation (or all-right spherical triangulation)
of the smooth manifold $(M,\cS)$. Let $\cS'$ be as in Theorem 1.1.1
Then, for every $\sigma\in K$ we have that $f(\sL(\sigma,K))$
is a smooth submanifold of $(M,\cS')$.}\\

\noindent {\bf Corollary 1.1.3.} {\it Let $M$, $\cS$ and $\cS'$ as in Theorem 1.1.1.
(or its addendum).
Then $K$ is $PL$-homeomorphic to $(M,\cS')$.}\\

\noindent {\bf Proof.} Since $K$ is a smooth
cubulation of $(M,\cS)$ we have $K\cong_{PL}(M,\cS)$.
On the other hand, by Theorem 1.1 we get $(M,\cS)\cong_{DIFF}(M,\cS')$.
Hence $K\cong_{PL}(M,\cS')$. This proves the corollary.\\

\noindent {\bf Remark.}
 Note that the image of the chart $h^\bullet_{\sigma}$
is the open normal neighborhood $\stackrel{\circ}{\s{N}}_\epsilon(\dsigma,K)$
of width $\epsilon$ of $\dsigma$ in $K$.
Even though we are assuming, for simplicity, that $\epsilon<1/2$ ($\epsilon<\pi/4$
in the spherical case) it can be checked from the proof of the Theorem 1.1.1 in \cite{O1} that we can
actually take $\epsilon=1$ ($\epsilon=\pi/2$) for the charts. 

\vspace{.4in}

\noindent {\bf 1.2. Induced Link Smoothings.}\\

Let $K$ be a cubical or all-right spherical complex. Then the links of $\sigma\in K$ are
all-right-spherical complexes. 
We explain here how to obtain from a given a collection of link smoothings for $K$ (and its corresponding normal atlas and structure) a collection of links smoothings for a link in $K$
(and its corresponding normal atlas and structure).\\

The all-right-spherical structure on 
$\sL(\sigma,K)$ induced by $K$ has all-right-spherical simplices
$\big\{ \,\tau\,\,\cap\,\, \sL(\sigma,K)\,\,  ,\,\, \tau\in K\big\} .$ 
Note that $\tau\,\,\cap\,\, \sL(\sigma,K)$ is non-empty only when $\sigma\subsetneq \tau$, hence we can write\\

\noindent\hspace{1.8in}$\sL(\sigma,K)=\big\{ \,\tau\,\,\cap\,\, \sL(\sigma,K)\,\,  ,\,\, \sigma\subsetneq \tau\in K\big\}. $\\

Since $\tau\,\,\cap\,\, \sL(\sigma,K)$ is a simplex in the 
all-right spherical complex $\sL(\sigma,K)$ we can consider its
link $\sL\,\Big(\, \tau\,\cap\,\sL\,(\sigma,K)  , \sL\,(\sigma,K)  \,\Big)$.
By definition we have:\\

\noindent \hspace{1.6in}$\sL\,\Big(\, \tau\,\cap\,\sL\,(\sigma,K)  , \sL\,(\sigma,K)  \,\Big)\,=\,\sL\,\big( \, \tau   , K  \,\big),$\\

\noindent provided we choose the radii and bases of the links
properly. In the formula above  radii and bases
are not specified but the radii are certainly not equal. 
The simple relationship between these radii is given by equation (1)
in the proof of Lemma 1.2 \cite{O1} (or the corresponding one in the spherical case; see Remark 1 after the proof of Lemma 1.2 \cite{O1}). 
For the cone link we have a similar formula but it is not
an equality, it is just an  identification, which we call $\Re$:\\

\noindent {\bf (1.2.1)}\hspace{1in}$\rC\sL \,\Big(\, \tau\,\cap\,\sL\,(\sigma,K)  , \sL\,(\sigma,K)  \,\Big)\,\xrightarrow{\,\,\,\,
\Re\,\,\,\,\,}\rC\sL\,\big( \, \tau   , K  \,\big).
$\\

\noindent On the term in the right side the radial segments are
``straight" (geodesic in each simplex of $K$) but on the left
side the radial segments are ``curved" (they lie in $\sL\,(\sigma,K)$).\\

%In particular we can write $\rC\sL\,\big( \, \tau   , K  \,\big)\sbs
% \sL\,(\sigma,K)$, when $\sigma\sbs\tau$.\\

\noindent {\bf Remark 1.2.2.} In the cubical case the identification (1.2.1) above can be done in the following way.
Let $v\in\rC\sL\,\big( \, \tau   , K  \,\big)=\rC\sL_r\,\big( \, \tau   , K  \,\big)$
at some $p\in \dot{\tau}$ and consider
$ \sL\,(\sigma,K)= \sL_s\,(\sigma,K)$ at some point $q\in \dsigma$, with
$d_K(p,q)=s$ with the segment $[p,q]$ perpendicular to $\sigma$.
Then $v$ corresponds to the point $v'$ in the segment $[q,v]$ at a distance
$s$ from $q$. The (angular) distance in $\sL(\sigma,K)$ from $p$ to $v'$ is
$tan^{-1}(\frac{d_K(v,p)}{s})$.
We shall
write $\Re:\rC\sL(\tau,K)\hookrightarrow \rC\sL \,\Big(\, \tau\,\cap\,\sL\,(\sigma,K)  , \sL\,(\sigma,K)  \,\Big)\sbs\sL(\sigma,K)$ for the radial
projection described above. And we will write $\Re=\Re\0{p,q,r,s}$ if we want to make explicit the
dependence of $\Re$ on the choices above.\\

\noindent {\bf Remark 1.2.3.} In the all-right spherical case
the identification $\Re$ is similar; we only need to replace
the formula $tan^{-1}(\frac{d_K(v,p)}{s})$ by
$tan^{-1}(\frac{tan\,d_K(v,p)}{sin\,s})$. The latter formula
is obtained using spherical trigonometry.\\

\noindent {\bf Remark 1.2.4.} Note that the representation of the radial projection $\Re$ in the chart
$\big(\, h_{\tau}^\bullet, \D^{n-j}\times\dot{\tau}\,\big)$, $\tau$ a $j$-simplex, is smooth.\\

Now assume in addition that $K$ is $PL$ manifold, and let
$\{ h_\sigma\}_{\sigma\in K}$ be a set of link smoothings for
(the links of) $K$.
Since we have $\sL\,\Big(\, \tau\,\cap\,\sL\,(\sigma,K)  , \sL\,(\sigma,K)  \,\Big)\,=\,\sL\,\big( \, \tau   , K  \,\big)$
we can say that the set of link smoothings $\{ h_\sigma\}_{\sigma\in K}$ for
$K$ induces, just by restriction, a set of link smoothings for $\sL(\sigma, K)$,
$\sigma\in K$, given by
$\{ h_\tau\}_{\sigma\subsetneq\tau}$. We have the following diagram:

$$\begin{array} {ccc}\bS^{n-i-1}&\xrightarrow{\,\,\,\,\,\,\,\,\,\,h\0{\tau^i\cap\sL\,(\sigma,K)}}
&\sL\,\Big(\, \tau\,\cap\,\sL\,(\sigma,K)  , \sL\,(\sigma,K)\Big) \\
\big|\big|&&\big|\big|\\\bS^{n-i-1}&\xrightarrow{\,\,\,\,\,\,\,\,\,\,\,\,\,\,\,h\0{\tau^i}\,\,\,\,\,\,\,\,\,\,\,\,\,\,\,}
&\sL\,\big( \, \tau   , K  \,\big)
\end{array}$$

\noindent The double vertical lines on the sides are equalities.
This diagram is a triviality in the sense that the function
on the top row $h\0{\tau^i\cap\sL\,(\sigma,K)}$ is equal to the function
on the bottom row $h\0{\tau^i}$; but we shall use the notation
$h\0{\tau^i\cap\sL\,(\sigma,K)}$ when we consider the link smoothing
for the link $\sL\,\big(\, \tau\,\cap\,\sL\,(\sigma,K)  , \sL\,(\sigma,K)\big)$,
instead of the link smoothing for the link $\sL\,\big( \, \tau   , K  \,\big)$.\\

Note that the atlas  $\cA_\sigma=\cA_{\sL(\sigma,K)}=\big\{h^\bullet_
{\tau\cap\sL\,(\sigma,K)}\big\}_{\sigma\subsetneq \tau}$ 
is a (a priori just topological) normal atlas on $\sL(\sigma,K)$.\\

We now change the notation slightly, to match the one we
will use in Section 1.3: we replace $\tau$ by $\sigma^i$. 
So, let $\sigma^k\sbs\sigma^i\in K$.
We mentioned above that we have $h\0{\sigma^i\cap\sL\,(\sigma^k,K)}=
h\0{\sigma^i}$; but because the map $\Re$ of (1.2.1)
is not an equality, we can not say that the neighborhood smoothings
$h\0{\sigma^i\cap\sL\,(\sigma,K)}^\bullet$ and
$h\0{\sigma^i}^\bullet$ are equal
(on their respective domains). But there is a relationship between $h\0{\sigma^i\cap\sL\,(\sigma^k,K)}^\bullet$ and $h\0{\sigma^k}$. The map $h\0{\sigma^i\cap\sL\,(\sigma^k,K)}^\bullet$ is given by the following composition:\\

\noindent {\bf (1.2.5)}\,\,\,\,\,\,\,
{\small $\D^{n-i} \times \big(\dsigma^i\cap\sL\,(\sigma^k,K)\,\big)\,\,
%\hookrightarrow
%\D^{n-i} \times \dsigma^i\,\,
\,\stackrel{{\mbox{{\tiny $ h^\bullet_{\sigma^i}$\,}}}}{\longrightarrow}\,\rC\sL\big(\sigma^i,K\big)\times\big(\dsigma^i\cap\sL\,(\sigma^k,K)\big)\,\,
\stackrel{{\mbox{{\tiny $\Re\times 1$}}}} {\longrightarrow}\, \sL(\sigma^k,K),\,$}\\

\noindent where $1$ is the identity on $\dsigma^i\cap\sL\,(\sigma^k,K)$,
$\Re$ is the map in Remark 1.2.2. Also the first arrow is really
the map $ h^\bullet_{\sigma^i}$ restricted to $\D^{n-i} \times (\dsigma^i\cap\sL\,(\sigma^k,K)\,)$ (the domain of $h^\bullet_{\sigma^i}$ is
$\D^{n-i} \times \dsigma^i$).
\vspace{.3in}

\noindent {\bf 1.3. Change of charts and smooth compatibility.} \\

In this section we define what a ``smoothly compatible
set of link smoothings" is and prove that a normal atlas is smooth
if and only if the set of link smoothings is smoothly compatible.\\

Let $f:K\ra M^n$ be a smooth cubulation (or all-right spherical triangulation).
Since in what follows of this subsection the function $f$ is not essential,
to simplify our notation we identify $K$ and $M$ via $f$. Let 
$\big\{h_{\sigma}\big\}$ be a set of link smoothings on $K$. Recall that this set determines a 
(not necessarily smooth) atlas $\cA=\Big\{h^\bullet_{\sigma}\Big\}$ on  $K$. \\

Let $\sigma^k\sbs\sigma^i\in K$. We say that the link smoothings
$h_{\sigma^k}$, $h_{\sigma^i}$ are {\it smoothly compatible} if
the neighborhood smoothing $h^\bullet_{\sigma^i\cap\sL\,(\sigma^k,K)}$
$$
h^\bullet_{\sigma^i\cap\sL\,(\sigma^k,K)}:\D^{n-i}\times\big(\dsigma^i
\cap\sL\,(\sigma^k,K)\big)
\xrightarrow{\,\,\,\,\,\,\,\,\,\,\,\,}\sL\,(\sigma^k,K)
$$
\noindent is a smooth embedding. Here we are considering 
$\sL\,(\sigma^k,K)$ with the smooth structure induced
by the link smoothing $h_{\sigma^k}:\bS^{n-k-1}\ra\sL(\sigma^k,K)$.
That is, we consider $\sL\,(\sigma^k,K)$ with the smooth structure
$(h_{\sigma^k})_*(\cS_{\bS^{n-k-1}})$, where $\cS_{\bS^{n-k-1}}$ is the canonical
smooth structure on $\bS^{n-k-1}$. Equivalently, 
$h_{\sigma^k}$, $h_{\sigma^i}$ are smoothly compatible if
the composition

$$
\D^{n-i}\times\big(\dsigma^i
\cap\sL\,(\sigma^k,K)\big)
\xrightarrow{\,\,\,\,\,h^\bullet_{\sigma^i\cap\sL\,(\sigma^k,K)}\,\,\,\,\,\,\,}\sL\,(\sigma^k,K)\xrightarrow{\,\,\,\,\,\,\,\,\,\,h_{\sigma^k}^{-1}\,\,\,\,\,\,\,\,\,\,}\bS^{n-k-1}
$$
\noindent is a smooth embedding.\\

\noindent {\bf Lemma 1.3.1.} {\it Fix $\sigma^k\in K$ and assume
$h_{\sigma^k}$ is smoothly compatible with $h_{\sigma^i}$,
for every $\sigma^i\supset\sigma^k$. Then the atlas
$\cA_{\sigma^k}=\cA_{\sL(\sigma^k,K)}=\big\{h^\bullet_
{\sigma^i\cap\sL\,(\sigma,K)}\big\}_{\sigma^k\subsetneq \sigma^i}$ 
is a smooth normal atlas on $\sL(\sigma^k,K)$. Moreover,
the link smoothing
$$h_{\sigma^k}:\bS^{n-k-1}\ra\Big(\sL(\sigma^k,K)\,\,,\,\,
\cS_{\sigma^k}\Big)$$
\noindent is a diffeomorphism. Here $\cS_{\sigma^k}$ is the normal
smooth structure induced by the normal atlas $\cA_{\sigma^k}$.}\\

\noindent {\bf Proof.} It follows from the fact that the maps
$h_{\sigma^k}^{-1}\circ h^\bullet_{\sigma^i\cap\sL\,(\sigma^k,K)}$
are smooth embeddings. This proves the lemma.\\

The set of link smoothings $\big\{h_{\sigma}\big\}$ is {\it smoothly
compatible} if $h_{\sigma^k}$, $h_{\sigma^i}$ are smoothly compatible
whenever $\sigma^k\sbs \sigma^i\in K$.\\

To simplify our notation write $S=\sL(\sigma^k,K)$, 
$\sigma^i_{S}=\sigma^i\cap S$, and 
$\dsigma^i_{S}=\dsigma^i\cap S$.
%Notice that if $\sL(\sigma^k,K)=\sL_t(\sigma^k,K)$ is based at $p\in \sigma^k$, then
%$\sigma_{S}^i$ is the set of points $x$ in $\sigma^i$ with $d\0{\sigma^i}(p,x)=t$ and $[p,x]$
%perpendicular to $\sigma^k$ at $p$. \\
%In section 2 we mentioned that for $\sigma^k\sbs\sigma^j$, using
%a radial projection, we  get the identification
%$\Re:\rC\sL(\sigma^j,K)\hookrightarrow \rC\sL \,\Big(\, \sigma^j_S  , \sL\,(\sigma^k,K)  \,\Big)\sbs\sL(\sigma^k,K)$. 
Now using (1.2.5) we can say that the set of link smoothings is smoothly compatible if
for every $\sigma^k\sbs\sigma^i$ the composition maps\\

\noindent {\bf (1.3.2)}\hspace{.8in}
{\small $\D^{n-i} \times \dsigma\0{S}^i\,\,
\,\stackrel{{\mbox{{\tiny $ h^\bullet_{\sigma^i}$\,}}}}{\longrightarrow}\,\rC\sL\big(\sigma^i,K\big)\times\dsigma^i_S\,
\stackrel{{\mbox{{\tiny $\Re\times 1$}}}} {\longrightarrow}\, \sL(\sigma^k,K)\,\stackrel{{\mbox{\tiny $(h_{\sigma^k} )^{-1}$}}}{\longrightarrow}\,\bS^{n-k-1}$}\\

\noindent are smooth embeddings. 
Notice that the image of the composition of the first two arrows is the normal neighborhood 
$\rC\sL \,\Big(\, \sigma^j_S , \sL\,(\sigma^k,K)  \,\Big)\times  \dsigma^j_S$ of $ \dsigma^j_S$ in $\sL(\sigma^k,K)$. Note that the map given in (1.3.2) can also be written as\\

\noindent {\bf (1.3.3)}\hspace{1.2in}
{\small $\D^{n-i} \times \dsigma\0{S}^i\,\,\stackrel{{\mbox{{\tiny $\Re'$}}}}{\longrightarrow}\, S'
\,\stackrel{{\mbox{{\tiny $ h^\bullet_{\sigma^i}$\,}}}}{\longrightarrow}\,S= \sL(\sigma^k,K)\,\stackrel{{\mbox{\tiny $(h_{\sigma^k} )^{-1}$}}}{\longrightarrow}\,\bS^{n-k-1}$}\\

\noindent where $S'=(h^\bullet_{\sigma^i})^{-1}(S)$ and $\Re'$ is the representation of $\Re\times 1_{\sigma^i_S}$ in the chart
$h^\bullet_{\sigma^i}$ (which is always smooth, see Remark 1.2.4).
Note that the inverse of the map given in (1.3.2) (or 1.3.3) is just the change of charts
$(h_{\sigma^i}^\bullet )^{-1}\circ h_{\sigma^k}^\bullet $ restricted to an open subset of $\bS^{n-k-1}$,
plus the ``straightening map"  $(\Re')^{-1}$.\\

\noindent {\bf Proposition 1.3.4.} {\it The set of link smoothings $\{h_\sigma\}$ is smoothly
compatible if and only if the atlas $\cA\big(\{h_\sigma\}\big)$ is smooth.}\\

\noindent {\bf Proof.} If the atlas $\cA$ is smooth then all chart maps $h^\bullet_{\sigma}$ are embeddings
(with respect to the smooth structure generated by $\cA$).
Therefore  the composition given in (1.3.2)  above
is smooth. Note the ``identification" $\Re$
is also smooth, for it is smooth in the chart corresponding to $\sigma^i$ (see Remarks 1.2.2, 1.2.4).
Hence $\{h_\sigma\}$ is smoothly compatible.
We prove the converse by induction on the codimension of the skeleta. \\

Write $k+j=n$.
Suppose the set of smoothings $\{h_\sigma\}$ is smoothly compatible.
Denote by $W_{\sigma^i}$ the image of $h_{\sigma^i}^\bullet$.
Assume that we have proved that the atlas  $\cA_{j}=\Big\{ \big(h^\bullet_{\sigma^i},
\D^{n-i}\times\dsigma^i\big)    \Big\}_{k<i}$ is smooth, and we want to prove that
the atlas  $\cA_{j+1}=\Big\{ \big(h^\bullet_{\sigma^i},
\D^{n-i}\times\dsigma^i\big)    \Big\}_{k\leq i}$ is smooth.
Note that $\cA_{j}$ is a smooth atlas on the complement $K-K_k$ of the $k$-skeleton $K_k$. The difference between the atlas $\cA_j$ and the atlas $\cA_{j+1}$
are the charts with maps $h^\bullet_{\sigma^k}$, for all link smoothings
$h_{\sigma^k}$ of $k$-cubes or ($k$-simplices) $\sigma^k$. We prove the proposition by
proving that the following maps are smooth embeddings\\

\noindent\hspace{1.8in}
$h^\bullet_{\sigma^k}|_{(\D^k\times\dsigma^k)-(\{0\}\times\dsigma^k)}\ra
(M-M_j,\cA_j).$\\

\noindent Fix $\sigma^k$. The open sets 
$U_{\sigma^i}=S\cap W_{\sigma^i}$, $\sigma^i\supset\sigma^k$, form an open cover of $S$. 
Note that $U_{\sigma^i}$ is a normal neighborhood of $\sigma\0{S}^i=\sigma^k\cap\sigma^i$
in $S$.
Write $V_{\sigma^i}=(h_{\sigma^k}^\bullet)^{-1}(\rC^+ U_{\sigma^i})$
(here the vertex $o$ of the open cone $\rC^+ U_{\sigma^i}$ is the center of $S$). Let $u\in S$. 
Take $i$ so that
$u\in U_{\sigma^i}$. \\

\noindent {\it Claim.  The map \,\,\,
$h_{\sigma^k}^\bullet|\0{V_{\sigma^i}\times\dsigma^k}: V_{\sigma^i}\times
\dsigma^k\ra W_{\sigma^i}\sbs (K-K_k,\cA_{j})$ is an embedding. }\\

\noindent {\it Proof of the Claim.} Since  $h^\bullet_{\sigma^i}$ is already a smooth embedding
it is enough to prove that the map $h=\big(h^\bullet_{\sigma^i}\big)^{-1}\circ
h_{\sigma^k}^\bullet|\0{V_{\sigma^i}} $ is a smooth embedding. But 
we can consider $V_{\sigma^i}\sbs \D^{n-k}-\{0\}=\bS^{n-k-1}\times(0,1)$. Write $\sigma^i=\sigma^l\times\sigma^k$, and let the link $S=\sL(\sigma^k,K)$ be based at $p\in\dsigma^k$. For $v=
(t\,u,w)\in V_{\sigma^i}\times \dsigma^k$, $u\in\bS^{n-k-1}$, we can write\\
 
\hspace{1.4in}$h(v)=( \alpha_u(t) ,w)\in (\D^{n-i}\times \dsigma^l)\times\dsigma^k=\D^{n-i}\times\dsigma^i,$\\

\noindent where $\alpha_u$ is the segment $[p,g(u)]$
and $g$ is the inverse of the map in (1.3.3). This proves the claim.\\

Since the open sets $V_{\dsigma^i}\times\dsigma^k$ cover 
$(\D^{n-k}\times\dsigma^k)-(\{0\}\times\dsigma^k)$ we can conclude that $h_{\sigma^k}^\bullet$ is a smooth embedding
away from $\{0\}\times\dsigma^k$, into the smooth manifold $(K-K_k,\cA_{j+1})$.
This proves the proposition.\\

\noindent {\bf Corollary 1.3.5.} {\it Let $\{h_\sigma\}$ be a set of
link smoothings on $K$, and let $\sigma^k\in K$.
If the atlas $\cA=\cA\big(\{h_\sigma\}_{\sigma\in K}\big)$ is smooth, then:
\begin{enumerate}

\item[(1)] The atlas
$\cA_{\sigma^k}=\cA_{\sL(\sigma^k,K)}=\big\{h^\bullet_
{\sigma^i\cap\sL\,(\sigma^k,K)}\big\}_{\sigma^k\subsetneq \sigma^i}$ 
is a smooth normal atlas on $\sL(\sigma^k,K)$. 
\item[(2)] The set of link smoothings $\{h_
{\sigma^i\cap\sL\,(\sigma,K)}\}_{\sigma^k\subsetneq \sigma^i}$ 
for the links of $\sL(\sigma^k,K)$ is smoothly compatible.

\item[(3)] The link smoothing
$$h_{\sigma^k}:\bS^{n-k-1}\ra\Big(\sL(\sigma^k,K)\,\,,\,\,
\cS_{\sigma^k}\Big)$$
\noindent is a diffeomorphism.

\item[(4)] Let $\cS'$ be the normal smooth structure
on $K$ induced by $\cA$, and let
$\cS_{\sigma^k}$ be the normal smooth structure
on $\sL(\sigma^k,K)$ induced by $\cA_{\sigma^k}$. We have that$$\cS'\big|_{\sL(\sigma^k,K)}\,\,=\,\,\cS_{\sigma^k}.
$$
\noindent Here $\cS'\big|_{\sL(\sigma^k,K)}$
denotes the restriction of $\cS'$ to $\sL(\sigma^k,K)$.
(Recall that, by Corollary 1.1.2 
$\sL(\sigma^k,K)$ is a smooth submanifold of
$(K,\cS')$.)

\end{enumerate}}

\noindent {\bf Proof.} Since $\cA$ is smooth, by Proposition 1.3.4,
the set of link smoothings $\{h_\sigma\}_{\sigma\in K}$
is smoothly compatible. This together with
Lemma 1.3.1 imply (1). Item (2) follows from 1(1) and Proposition 1.3.4
(applied to the complex $\sL(\sigma^k,K)$).
Item (3) also follows from Lemma 1.3.1 and (2)
(i.e. the fact that $\{h_
{\sigma^i\cap\sL\,(\sigma,K)}\}_{\sigma^k\subsetneq \sigma^i}$  is smoothly compatible). Item (4) follows from (3)
and the fact that
$$h_{\sigma^k}:\bS^{n-k-1}\ra\big(K\,\,,\,\,
\cS'\big)$$
\noindent is a smooth embedding (see Remark 1 before Theorem 1.1.1). This proves the corollary\\

%And, since a link smoothing $h_\sigma$ is the restriction of the embedding $h^\bullet_\sigma$ we get the following corollary.\vspace{.1in} 

%\noindent {\bf Corollary 2.6.} {\it For every $\sigma^k\in K^n$, the link smoothing $h_{\sigma}:\bS^{n-k-1}\ra \big(\sL(\sigma,K),\cS'_\sigma\big)$ is a diffeomorphism.}
\vspace{.3in}

\noindent {\bf 1.4. A few technical results.}\\

This is a technical subsection. We present some results that will be needed later.\\

Let $f:K\ra M$ be a smooth cubulation  of $M$ and $\cA$ be a normal atlas for $K$, inducing the smooth structure $\cS'$ on $M$. In general for a (closed)
cube (or simplex) $\sigma$ the inclusion $\sigma\hookrightarrow (M,\cS')$ is (almost always) not
a smooth embedding (see \cite{O1}). But we prove in the
next lemma that a weaker regularity condition holds. Consider $\sigma^j=\sigma^l\times\sigma^k\in K$.  As before we identify a normal
neighborhood of $\sigma^k$ in $\sigma^j$ with 
$\rC\sL(\sigma^k,\sigma^j)\times \sigma^k$ and write an element
in $\rC\sL(\sigma^k,\sigma^j)$ in the form $tu$, $u\in\sL(\sigma^k,\sigma^j)$.
We have an inclusion $\sL(\sigma^k,\sigma^j)\times \sigma^k\sbs
\rC\sL(\sigma^k,\sigma^j)\times \sigma^k$.
Also denote the inclusion $\sigma^j\hookrightarrow (M,\cS')$ by $\iota$.\\

\noindent {\bf Lemma 1.4.1.} {\it 
Let $(u\0{0},p\0{0}), (u_n, p_n)\in \sL(\sigma^k,\sigma^j)\times \dsigma^k\sbs
\sigma^j$ with $(u_n,p_n)\ra (u\0{0},p\0{0})$.
Also let $t_n\ra 0\in [0,\infty)$, and let $(0, v\0{0}), (0,v_n)\in T_{(u\0{n},p\0{n})}\Big(\rC\sL(\sigma^k,\sigma^j)\times \sigma^k\Big)$ (hence
they are parallel to $\sigma^k$), with $v_n\ra v\0{0}$. Then}
\begin{enumerate}
\item[{\it (i)}] {\it  we have that $D\iota\0{(t_nu_n,p_n)}(u_n,0)\ra D\iota\0{(0,p\0{0})}(u\0{0},0)$,}
\item[{\it (ii)}] {\it  we have that $D\iota\0{(t_nu_n,p_n)}(0,v_n)\ra D\iota\0{(0,p\0{0})}(0,v\0{0})$.}
\end{enumerate}

\noindent {\bf Proof.} Just take the chart $h^\bullet_{\sigma^j}$ and recall that
$h^\bullet_{\sigma^j}$ is a product map and respects the radial structure. This proves the lemma.\\

Let $U$ be a bounded open set of $\R^n$. A smooth map $f:U\ra\R^N$
is {\it polynomially bounded with respect to a subset } $B\sbs \R^n$ if for every
$p\in B$ and every partial derivative $\p^\alpha f$ of $f$
we have constants $C,\, m\in \R$ such that $|\p^\alpha f(p)|\leq C d\0{\R^n}(p, B)^m$.
The map is {\it $C^1$ well-bounded} if the norm of the first derivative $|Df|$ is bounded
and bounded away from zero.
For $U\sbs\bS^n$ we write $\rC^+U=\rC U-\{0\}$,
where $\rC U$ is the cone $\rC U=U\times [0,1]/U\times \{0\}$
of $U$.\\

\noindent {\bf Lemma 1.4.2.} {\it Let $U\sbs\bS^n$ be open and $f:U\ra f(U)\sbs\R^N$ be smooth. Let $V$ open with $\bar{V}\sbs U$. Then the cone map
$\rC \big(f|_{V}\big):\rC^+V\ra \rC f(V)\sbs \R^N$ is polynomially bounded at 0.
Furthermore, $\rC \big(f|_{V}\big)$ is also $C^1$ well-bounded.}\\

\noindent {\bf Proof.} We have $(\rC f) (x)=|x|f(\frac{x}{|x|})$, and the lemma follows from differentiating this equation.
This proves the lemma.\\

\noindent {\bf Remark.} Actually since $\rC f$ is a cone map then $D(\rC f)_x=D(\rC f)\0{\frac{x}{|X|}}$, for $x\neq 0$.
In particular
if $u\in V$, and $t> 0$ then
$D(\rC f)_{tu}u=f(u)$.\\

Let  $\cA=\Big\{ \big(  h_{\sigma^i}^\bullet,\D^{n-i}\times\dsigma^i    \big)\Big\}$
be a normal atlas. Let $S=f(\sL(\sigma^k,K))$ and $U'\sbs S\cap W_{\sigma^i}$,
where $\sigma^i>\sigma^k$ and $W_{\sigma^i}$ is the image of 
$h_{\sigma^i}^\bullet$. Write $U=\big(h^\bullet_{\sigma^k}\big)^{-1}(U')$
and let $V$ open with $\bar{V}\sbs U$. We have the following corollary about
the change of charts $\big(h^\bullet_{\sigma^i}\big)^{-1}\circ h_{\sigma^k}^\bullet:
\rC^+V\ra \D^{n-i}\times\sigma^i\sbs\R^n$, restricted to the cone of $V$.\\

\noindent {\bf Corollary 1.4.3.} {\it The change of charts
$\big(h^\bullet_{\sigma^i}\big)^{-1}\circ h_{\sigma^k}^\bullet$,
restricted to $\rC^+ V$, is polynomially bounded at 0, and
$C^1$ well-bounded.}

\vspace{.5in}

\noindent {\bf 1.5. The case of manifolds with codimension zero singularities.}\\

Here we treat the case of manifolds with a one point singularity.
The case of manifolds with many (isolated) point singularities
is similar. \\

Let $Q$ be a smooth manifold with a one point singularity $q$, that is
$Q-\{q\}$ is a smooth manifold and there is a topological embedding
$\rC_1 N\ra Q$, with $o\0{\rC N}\mapsto q$, that is a smooth embedding outside 
the vertex $o\0{\rC N}$. Here $N=(N,\cS_N)$ is a closed smooth manifold (with smooth
structure $\cS_N$).
Also $\rC_1 N $ is the (closed) cone of width 1 and we identify
$\rC_1 N-\{o\0{\rC N}\}$ with $N\times (0,1]$. We write $\rC_1N\sbs Q$.
We say that the {\it singularity $q$ of $Q$ is  modeled on  $\rC N$.}\\

Assume $(K,f )$ is a smooth cubulation of $Q$, that is
\begin{enumerate}
\item[(i)] 
$K$ is a cubical complex.
\item[(ii)] $f:K\ra Q$ is a homeomorphism. Write $f(p)=q$ and $L=\sL(p,K)$.
\item[(iii)] $f|_{\sigma}$ is a smooth embedding for every
cube $\sigma$ not containing $p$.
\item[(iv)]  $f|_{\sigma-\{p\}}$ is a smooth embedding for every
cube $\sigma$ containing $p$.
\item[(v)]  %There is $g:L\ra N$ such that $(L,g)$ is a smooth triangulation of $N$. That is, 
$L$ is $PL$ homeomorphic to $(N,\cS_N)$.
\end{enumerate}

Many of the definitions and results given before
for smooth cube manifolds still hold (with minor changes)
in the case of manifolds with a one point singularity:\\

\begin{enumerate}
\item[{\bf (1)}] A {\it link smoothing} for $L=\sL(p,K)$ (or $p$) is just a
homeomorphism $h_p:N\ra L$.
Since all but one of the links of $K$ are spheres, sets of  link smoothings for $K$
are defined, that is they are sets of link smoothings for the sphere links plus
a link smoothing for $L$.
\item[{\bf (2)}] Given a set of link smoothings for $K$ we get a set of normal charts
as before. For the vertex $p$ we mean the cone map 
$h_p^\bullet=f\circ\rC h_p:\rC N\ra Q$. We will also denote the restriction
of $h_p^\bullet$ to $\rC N-\{o\0{\rC N}\}$ by the same notation $h_p^\bullet$.
As before $\{h^\bullet_\sigma\}_{\sigma\in K}$ is a {\it (topological) normal atlas on $Q$
with respect to $K$}. The atlas on $Q$ is {\it smooth} if all transition functions are smooth, where for the case $h_p^\bullet:\rC N-\{o\0{\rC N}\}\ra Q-\{q\}$ we are identifying $\rC N-\{
o\0{\rC N}\}$ with $N\times (0,1]$ with the product smooth structure obtained from
\s{some} smooth structure ${\tilde{\cS}}_N$ on $N$.
A smooth normal atlas on $Q$ with respect to $K$ induces, by restriction,  a smooth normal structure on $Q-\{q\}$ with respect to $K-\{p\}$ (this makes sense even though
$K-\{p\}$ is not, strictly speaking, a cube complex).
\item[{\bf (3)}] We say that the set $\{h_\sigma\}$ is smoothly compatible if
condition (1.3.2) holds.  It 
is straightforward to verify that Proposition 1.3.4 holds:  $\{h_\sigma\}$ is smoothly compatible if and only if
$\{h^\bullet_\sigma\}$ is a smooth atlas on $K$. In this case we say that
the smooth atlas  $\{h^\bullet_\sigma\}$ (or the induced smooth structure, or the set $\{h_\sigma\}$) 
is {\it correct with respect to $N$} if $\cS_N$ and ${\tilde{\cS}}_N$ are diffeomorphic.
\item[{\bf (4)}] Also it is straightforward to verify that
Corollary 1.3.5 holds in our present case. 
\item[{\bf (5)}] Theorem 1.1.1 also holds in this context.
\end{enumerate}

\noindent {\bf Theorem 1.5.1.} {\it  Let $Q$ be a smooth manifold
with one point singularity $q$ modeled on $\rC N$, where $N$ is a closed
smooth manifold. Let $(K,f)$ be a smooth cubulation of $Q$. Then
$Q$ admits a normal smooth structure with respect to $K$,  which restricted to $Q-\{q\}$ is diffeomorphic to $Q-\{q\}$.
Moreover, this normal smooth structure is correct with respect to $N$ if}
\begin{enumerate}
\item[(a)] {\it $dim\, N\leq 4$.}
\item[(b)] {\it $dim\, N\geq 5$ and the Whitehead group $Wh(N)$ of $N$ vanishes. }
\end{enumerate}

This Theorem is proved in \cite{O1}.\vspace{.6in}

\begin{center} {\bf \large 2. Charney-Davis Strict Hyperbolization Pieces.}
\end{center}

We use some of the notation in \cite{ChD}. In particular the canonical $n$-cube
$[0,1]^n$ will be denoted by $\square^n$.
(This differs with the notation used in Section 1, where an $n$-cube was denoted by
$\sigma^n$.)
Also $B_n$ is the isometry group
of $\square^n$.\\

A {\it Charney-Davis strict hyperbolization piece of dimension $n$} is a compact
connected orientable hyperbolic $n$-manifold with corners satisfying the properties stated in Lemma 6.2 of \cite{ChD}.
The group $B_n$ acts by isometries on $X^n$
and there is a smooth map $f:X^n\ra\square^n$ constructed in Section 5 of
\cite{ChD} with certain properties. 
We collect some facts from \cite{ChD}.

\begin{enumerate}
\item[(1)] For any $k$-face $\square^k$ of $\square^n$ we have that $f^{-1}
(\square^k)$ is totally geodesic in $X^n$ and it is a Charney-Davis hyperbolization
piece of dimension $k$. The totally geodesic submanifold (with corners)
$f^{-1}(\square^k)$ is a $k$-$face$ of $X^n$. Note that the intersection of
faces is a face and every $k$-face is the intersection of exactly $n-k$ distinct
$(n-1)$-faces.
\item[(2)] The map $f$ is $B_n$-equivariant.
\item[(3)] The faces of $X^n$ intersect orthogonally.
\item[(4)] The map $f$ is transversal to the $k$-faces of $\square^n$, $k<n$.
%\item[(5)]
\end{enumerate}

\noindent The $k$-face $f^{-1}(\square^k)$ of $X$ will be 
denoted by $X_{\square^k}$. The interior $f^{-1}(\dot{\square^k})$ will be denoted by $\dX_{\square^k}$.\\

\noindent {\bf Lemma 2.1.} {\it For every $n$ and $r>0$ there is a 
Charney-Davis hyperbolization piece of dimension $n$ such that the 
widths of the normal neighborhoods of every $k$-face, $k=0,...n-1$,
are larger than $r$.}\\

\noindent {\bf Proof.} A piece $X^n$ is constructed in Section 6 of 
\cite{ChD} by cutting a closed hyperbolic $n$-manifold $M$ along a
system $\{Y_i\}_{i=1}^n$ of codimension one totally geodesic submanifolds of $M$
that intersect orthogonally. The $Y_i$'s are
orientable and two sided.
The group $B_n$ acts by isometries on $M$, permuting the $Y_i$'s.
In particular each $Y_i$ is contained in the fixed point set of a nontrivial
 isometric involution $r_i$. Therefore $r_i$ interchanges
both sides of $Y_i$.\\

The $(n-1)$-faces of $X^n$ correspond to the $Y_i$'s and a $k$-face of
$X^n$ corresponds to the (transverse) intersection of $n-k$ different $Y_i$'s. \\

\noindent {\it Claim.} {\it It is
enough to show that the $Y_i$'s have normal neighborhoods with large width.}\\

\noindent {\it Proof of Claim.} Let $Z=Y_i\cap Y_j$, $Y_i\neq Y_j$, and assume both $Y_i$, $Y_j$
have normal neighborhoods of width larger than $r$.
Let $\alpha$ be a path with end points in $Z$ of length $<2r$. Then $\alpha$
lies in the normal neighborhood of $Y_i$.
Using the distance decreasing 
normal geodesic deformation of $Y_i$ we can deform $\alpha$, rel end points,
to a shorter path $\beta$ in $Y_i$. Repeat the same argument now with
$\beta$, and using the fact that $Y_i$ and $Y_j$ intersect orthogonally,
we get that the deformation of $\beta$ lies in $Z$. 
The proof for larger intersections is similar. This proves the claim.\\

We continue with the proof of Lemma 2.1. The claim above and the existence
of the nontrivial isometric involutions $r_i$ imply that it is enough to have
$M$ with large injectivity radius. (To see this let $\alpha:[0,1]\ra
M$, $\alpha(0),\alpha(1)\in Y_i$, of length $<2r$, with $\alpha$ not homotopic,
$rel\,\{0,1\}$, to a path in $Y_i$. Then take $\beta=
(r_i\circ\alpha)^{-1}*\alpha$, which is a non-nullhomotopic loop
of length $<4r$.)\\

To obtain $M$ with large injectivity radius recall that $M$ is given in
Section 6 of \cite{ChD} as $M=\HH^n/\Gamma$, where
$\Gamma=\Gamma({\cal{J}})$ is a congruence subgroup
given by the ideal ${\cal{J}}$ of the ring of integers of the totally real
quadratic extension $K=\Q(\sqrt{d})$ of $\Q$. The only conditions required
for ${\cal{J}}$ are that $\Gamma$ is torsion free and $\Gamma\sbs SO_o(n,1)$.
Hence any deeper congruence subgroup $\Gamma'$ will serve as well.
But, by taking deeper congruence subgroups we can increase (in a well-known way)
the injectivity radius as much as we want: let $\gamma_1,...,\gamma_k\in\Gamma$ correspond 
to the closed geodesics of length $<r$, and take a deeper 
ideal that contains none of the $\gamma_i$'s. This proves the lemma.\\

Let ${\square^k}\sbs {\square^n}$.
We will denote the normal neighborhood of $X_{\square^k}$
in $X_{\square^n}$ by $\s{N}\0{s}(\dX_{\square^k})=\s{N}\0{s}(\dX_{\square^k},\dX_{\square^n})$. That is, it is
the union of 
all geodesics of length $\leq s$ in the hyperbolization piece
$X_{\square^n}$ that begin at, and are normal to, $\dX_{\square^k}$. In the next two results we 
assume that the widths of the normal neighborhoods of all 
faces $X_{\square^k}$ of $X_{\square^n}$ are
greater than some large $ s\0{0}$. \vspace{.1in}

\noindent {\bf Lemma 2.2.} {\it 
We can find hyperbolization pieces $X_{\square^n}$ such that the following holds. Let 
 $\square^k=\square^i\cap\square^j$,  $k\geq 0$.
Let $s\0{1}, s\0{2}>0$ with $s\0{1}+ s\0{2} < s\0{0}$. Then}\,\,\, $
\sN\0{s\0{1}}\big(\dX_{\square^i}  \big)\,\, \cap\,\,\sN\0{s\0{2}}\big(\dX_{\square^j}  \big)
\,\,\, \sbs\,\,\, \sN\0{s\0{1}+ s\0{2}}\big(\dX_{\square^k}  \big)
$.\vspace{.1in}

\noindent {\bf Proof.} We prove the result for
two $(n-1)$-faces $\square^k=\square^{n-1}_1$ and $\square^j=\square^{n-1}_2$.
The general result has a similar proof. We use the construction and notation given in Lemma 2.1. Therefore
it is enough to prove that
$\sN\0{s\0{1}}\big(Y_1  \big)\,\, \cap\,\,\sN\0{s\0{2}}\big(Y_2  \big)
\,\,\, \sbs\,\,\, \sN\0{s\0{1}+ s\0{2}}\big(Y_1\cap Y_2  \big)
$, were $\sN\0{s\0{1}}(Y )$ is the $s$-normal
neighborhood of $Y$ in $M$. We assume that the injectivity radius of $M$ is $>4s\0{0}$.
Let $p\in \sN\0{s\0{1}}\big(Y_1  \big)\,\, \cap\,\,\sN\0{s\0{2}}\big(Y_2  \big)$. Therefore there are  length minimizing geodesic 
segments $[u_i,p]$ in $M$ from $u_i\in Y_i$ to $p$, perpendicular to
$Y_i$, and of lengths $\leq s\0{i}$, $i=1,2$. Take liftings
$p',\, u_i', [u_i',p']$ of $p,\, u_i, [u_i,p]$, respectively,
to the universal cover $\HH^n$ of $M$.
Also let $Y_i'$ be the component of the pre-image 
of $Y_i$ that contains $u_i'$. We claim that
$Y_1'$ and $Y_2'$ interesect (othogonally)
in a non-empty subset. Indeed, if
$Y_1'\cap Y_2'=\emptyset$ there is a
geodesic segment $[p_1',p_2']$, $p_i'\in Y_i'$, of length $\leq s\0{1}+s\0{2}$
orthogonal to both $Y_i$'s
(because $[u_1',p']\ast[p',u_2']$ joins $Y_1'$ to $Y_2'$ and has
length $\leq s\0{1}+s\0{2}$). Then the image 
$[p_1,p_2]$  of
$[p_1',p_2']$ in $M$ is orthogonal to both $Y_i$'s.
As in the proof of Lemma 2.1, using the isometric
involutions $r_i$ and the segment
$[p_1,p_2]$ we obtain a non-trivial closed geodesic
of length $\leq 4(s\0{1}+s\0{2})<4s\0{0}$, which
is a contradiction because the
injectivity radius of $M$ is $>4s\0{0}$.
Therefore $Y_1'$ and $Y_2'$ have non-empty
intersection. \vspace{.1in}

Since $Y_1'$ and $Y_2'$ intersect
orthogonally, hyperbolic geometry arguments
show that 
there is a
geodesic segment $[u_1',q]$ in $Y_1'$, 
$q\in Y_1'\cap Y_2'$, of length $\leq s\0{2}$.
Since the length of the path
$[p',u_1']\ast[u_1',q]$ is $\leq s\0{1}+s\0{2}$
it follows that the distance from $p'$ to
$Y_1'\cap Y_2'$ is $\leq s\0{1}+s\0{2}$.
This proves the lemma.
 \vspace{.1in}

The next lemma gives a better estimate.\vspace{.1in}

\noindent {\bf Lemma 2.3} {\it Let 
 $\square^k=\square^i\cap\square^j$,  $k\geq 0$.
Let $s\0{1}, s\0{2}, s < \frac{1}{2}s\0{0}$ be positive real numbers such that
$\frac{sinh\, s\0{1}}{sinh\, s},\frac{sinh\, s\0{2}}{sinh\, s} \leq\frac{\sqrt{2}}{2}$. Then}\,\,\, $
\sN\0{s\0{1}}\big(\dX_{\square^i}  \big)\,\, \cap\,\,\sN\0{s\0{2}}\big(\dX_{\square^j}  \big)
\,\,\, \sbs\,\,\, \sN\0{s}\big(\dX_{\square^k}  \big)
$.\vspace{.1in}

\noindent {\bf Proof.} We use the objects and notation from
the proof of Lemma 2.2.
At the end of the proof of Lemma 2.2 we found a
geodesic segment $[u_1',q]$ in $Y_1'$, 
$q\in Y_1'\cap Y_2'$, and perpendicular to $ Y_1'\cap Y_2'$.
Similarly, we can find a
geodesic segment $[u_2',q']$ in $Y_2'$, 
$q\in Y_1'\cap Y_2'$, and perpendicular to $ Y_1'\cap Y_2'$. 
Since we are in $\HH^n$ it is straightforward to show
that $q'=q$. 
Let $a_i$ be the length of $[u_i',p]$. Then $a_i\leq s\0{i}$.
Denote by $t$ the distance between $p$ and $q$.  We get right hyperbolic triangles with vertices
$p'$, $q$, $u_i'$ (right at $u_i'$), and hypotenuse equal to $t$. 
Let $\theta_i$ be the angle at $q$.
Thus $\theta_i$ is opposite to the side with length $a_i$.
By the hyperbolic law of sines we have $sin\, \theta_i=\frac{sinh\,a_i}{sinh\, t}$,
and by hypothesis we get\vspace{.1in}

\hspace{1.4in}{\small $
sin\,\theta_i\,\,=\,\, \frac{sinh\,a_i}{sinh\, t}\,\,\leq\,\,
\frac{sinh\,s_i}{sinh\, t}\,\,=\,\, \frac{sinh\,s_i}{sinh\, s}\,
\frac{sinh\, s}{sinh\, t}    \,\,\leq \,\,\frac{\sqrt{2}}{2}\,\frac{sinh\,s}{sinh\, t}
$}\vspace{.1in}

We want to prove that $t\leq s$. Suppose $t>s$. It follows then from the inequality
above that $sin\,\theta_i<\sqrt{2}/2$, thus $\theta_i<\pi/4$.
Let $S$ be the link of $Y_1'\cap Y_2'$ at $q$ (which is isometric 
to  $\bS^{n-k-1}$). The segments
$[q,u_i']$ intersect $S$ in two different vertices $v_{i}$. Since the sets $S\cap Y_i'$ are disjoint,  the
(angle) distance $d\0{S}(v_{1},v_{2})$ between $v_{1}$ and $v_{2}$
is a least $\pi/2$. Also the segment
$[p',q]$ intersects $S$ in a point $w$, and we have $\theta_{i}=d\0{S}(u,v_{i})$.
Consequently\vspace{.05in}

\hspace{1in}{\small $
\frac{\pi}{2}\,\,\leq\,\,  d\0{S}(v_{1},v_{2})\,\,\leq\,\, 
d\0{S}(v_{1},w)\,+\, d\0{S}(w,v_{2})\,\, =\,\, \theta\0{1}\,+\,\theta\0{2}\,\,
< \frac{\pi}{4}\,+\,\frac{\pi}{4}\,\,=\,\, \frac{\pi}{2}$} \vspace{.05in}

\noindent which is a contradiction. This proves the lemma.\vspace{.1in}
\vspace{.1in}

We need some extra properties for the map $f$, so we give an explicit construction of it.
Recall from the proof of Lemma 2.1 that $X$ is obtained from the closed hyperbolic manifold
$M$ by cutting along the system $\{Y_i\}$. Similarly the map $f$ is obtained in \cite{ChD}
from a map $\varphi :M\ra  \T^n$. And this map has coordinate maps
$\varphi_i:M\ra\bS^1$, which are constructed by applying the Pontryagin-Thom
construction to the framed (two sided) codimension one submanifolds
$Y_i$. Here we need a bit more details so we give a specific construction for
$\varphi$.\\

Let $Y_i\times (-r,r)\sbs M$ be the normal geodesic neighborhood of
$Y_i$ of width $r>2$. Hence for $p=(y,t)\in Y_i\times(-r,r)$, the smooth map $p\mapsto t(p)=t$ gives the signed distance to $Y_i$.
Let  $\eta:\R \ra [-1,1]$ be a non-decreasing smooth map such that $\eta (t)=t/r$ for $t\in (-r+1,r-1)$,\,
$\eta (t)=1$ for $t\geq r$,\, $\eta (t)=-1$ for $t\leq -r$. By identifying $(\bS^1,1)$ with
$\big([-1,1]\big/\{{\mbox{{\small $-1\sim 1$}}}\}\,\, ,\,\, 0\big)$,
the smooth map $\eta\circ t$ induces
the smooth map $\varphi_i:Y_i\times (-r,r)\ra \bS^1$ that can be
extended to the whole of $M$.  Note that $\varphi_i^{-1}(1)=Y_i$ and
$\varphi^{-1}(\T^{i-1}\times\{1\}\times\T^{n-i})=Y_i$.
After cutting along the $Y_i$'s we get the map $f:X\ra\square^n$ and each $Y_i$
corresponds to two $(n-1)$-faces $X_{\square\0{i,0}^{n-1}}$, $X_{\square\0{i,1}^{n-1}}$
(one for each side of $Y_i$), where $\square\0{i,j}^{n-1}=\square^{n}\cap
\{x_i=j\}$, $j=0,1$. Moreover, the normal neighborhood
$Y_i\times(-r,r)$ corresponds to the two one-sided normal neighborhoods
$X_{\square\0{i,0}^{n-1}}\times [0,r)$, $X_{\square\0{i,1}^{n-1}}\times [0,r)$.
Write $f(p)=(f_1(p),...,f_n(p))\in\square^n\sbs\R^n$. Then if $p\in 
X_{\square\0{i,0}^{n-1}}\times [0,r)$, we have $f_i(p)=\frac{1}{2}\eta (t_i(p))$, where $t_i(p)$ is the distance
to $X_{\square\0{i,0}^{n-1}}$. Similarly if $p\in X_{\square\0{i,1}^{n-1}}\times [0,r)$
we have $f_i(p)=1-\frac{1}{2}\eta (t_i(p))$, where $t_i(p)$ is the distance
to $X_{\square\0{i,1}^{n-1}}$. 
And if $p\in X_{\square\0{i,0}^{n-1}} \times [0,r-1)$, we have 
$f_i(p)=\frac{1}{2r}\,t_i(p)$, 
and similarly for $X_{\square\0{i,1}^{n-1}}$.
In particular $p\in X_{\square\0{i,0}^{n-1}} \times [0,a)$, $a\leq r-1$, if and only
if $f(p)\in\square\0{i,0}^{n-1} \times [0,\frac{a}{r})$.
In what follows of this paper we assume
$\varphi$ and $f$ are constructed as above.\\

In what follows we will write $\square\0{i}^{n-1}=\square\0{i,0}^{n-1}$, if the context is clear.\\

\noindent {\bf Proposition 2.4.} {\it The derivative of $f$ sends normal vectors
to $X_{\square^k}$ to normal vectors to $\square^k$.}\\

\noindent {\bf Proof.} For simplicity we assume that $k=n-2$. The proof for general
$k$ is similar. We can also assume that $\square^{n-2}=\square_1^{n-1}\cap\square_2^{n-2}$,
where $\square^{n-1}_i=\square^n\cap\{x_i=0\}$.
Write $U_{i,j}=X_{\square_{i,j}^{n-1}}$, $U_i=U_{i,0}$, and $W=X_{\square^{n-2}}=U_1\cap U_2$.
Let $p\in W$. We certainly have that $p\in U_1\times[0,r)$ and 
 $p\in U_2\times[0,r)$. 
 %Assume for simplicity that 
%$p\in U_3\times[0,r)$ but $p\notin U_i\times[0,r)$, $i\geq 4$.
%Let $u\in T_p W$ be normal to $W$ and $\alpha$ be the geodesic with
%$\alpha (0)=p$ and $\alpha'(0)=u$.
We have to prove that $(Df_i)_p\,u=0$ for $i\geq 3$, where $u$ is orthogonal to $X_{\square^{n-2}}$.
For each $i\geq 3$ we have two cases. Case 1:
$p\notin U_{i,j}\times [0,r)$ for $j=0$ and $j=1$. In this case
it follows that  $(Dt_i)_p=0$, hence  $(Df_i)_p=0$.
Case 2: $p\in U_{i,j}\times [0,r)$ for $j=0$ or $j=1$.
Say $p\in U_3\times[0,r)$.
We want to prove that $(Df_3)_p\,u=0$. Since $p\in U_3\times[0,r)$
there is a geodesic $\beta$ in $U_3\times[0,r)$, beginning at $U_3$, normal to $U_3$ 
and ending in $p$. Also $t_3(p)$ is equal to the length of $\beta$. Since
$\alpha$ and $\beta$ are perpendicular at $p$ the function $t\mapsto t_3(\alpha(t))$
has a minimum at 0. Hence $D(t_3)_p\, u=0$, therefore $D(f_3)_p \,u= \eta'(t_3(p))\,
D(t_3)_p\, u=0$. This proves the proposition.\\
\\

For a $k$-face $X_{\square^k}$ and $p\in X_{\square^k}$, the set of inward normal vectors to
$X_{\square^k}$ at $p$ can be identified with the canonical all-right $(n-k-1)$-simplex
$\Delta^{n-k-1}$. In this sense we consider $\Delta^{n-k-1}\sbs T_pX$.
Similarly we can consider $\Delta^{n-k-1}\sbs T_q\square^n$, for 
$q\in \square^k$. We make the convention that
the two identifications above are done with respect to an ordering of
the $(n-1)$-faces $X_{\square^{n-1}}$ of $X$ and the corresponding
ordering for $\square^n$. For instance the vectors in
$\Delta^{n-k-1}\sbs T_pX$ tangent to some $X_{\square^{n-1}}$
correspond to the same $(n-1)$-face of $\Delta^{n-k-1}$
as the vectors in $\Delta^{n-k-1}\sbs T_{f(p)}\square^n$ tangent to  
$\square^{n-1}$. With these identifications we get coordinates on
$\Delta^{n-k-1}$: we write $(u_1,...,u_n)=u\in \Delta^{n-k-1}$
where $u_i$ is the angle between $u$ and $X_{\square_i^{n-1}}$ (or $\square_i^{n-1}$).\\

\noindent {\bf Proposition 2.5.} {\it For $p\in \dX_{\square^k}$,
we have that
\begin{enumerate}
\item[(i)] $Df_p$ sends non-zero normal vectors to non-zero normal vectors.
\item[(ii)] For $u\in \Delta^{n-k-1}$ we have $Df_p(u)=\frac{1}{2r}u$.
\item[(iii)] {\bf n} $\circ \,(Df_p|\0{\Delta^{n-k-1}}):\Delta^{n-k-1}
\ra\Delta^{n-k-1}$
is the identity, where {\bf n}$(x)=\frac{x}{|x|}$ is the normalization map.
\end{enumerate}}

\noindent {\bf Remark.} In (iii) we are using  the coordinates on 
$\Delta^{n-k-1}$ mentioned above to identify the normal tangent spaces of
$X_{\square^k}$ and $\square^k$.\\

\noindent {\bf Proof.} 
For simplicity we assume that $k=n-3$. The proof for general
$k$ is similar. As in the proof of Proposition 2.4 write
$\square^{n-1}_i=\square^n\cap\{x_i=0\}$ and $U_i=X_{\square^{n-1}_i}$. For simplicity take
$W=X_{\square^{n-3}}=U_1\cap U_2\cap U_3$.
Let $u=(u\0{1},u\0{2},u\0{3})\in \Delta^{2}\sbs T_p\, X$.
%hence $u_i=0$, $i\geq 4$. 
Then $u_i$ is the spherical distance from $u$ to $T_p U_i$, $i=1,2,3$
(or angle between $u$ and $U_i$).
Let $\alpha$ be the geodesic with $\alpha(0)=p$ and $\alpha'(0)=u$. Then
$u_i$ is the angle, at $p$, between $\alpha$ and $U_i$.
%We have to prove that $Df_p u=u$. By proposition 2.4 $D(f_i)_p\, u=0$ %$i\geq 4$. 
We have to prove that $D(f_i)_p\, u=\frac{1}{2r}u_i$ for $i=1,2,3$. Since the length of $\alpha|_{[0,t]}$
is $t$ and the distance  from $\alpha (t)$ to $U_i$, $i=1,2,3$, is $t_i(\alpha(t))$,
we get a  right hyperbolic triangle with hypotenuse of length $t$
and side equal to $t_i(\alpha(t))$ with opposite angle $u_i$. Hence,
the hyperbolic law of sines implies

\begin{equation*}t_i(\alpha(t))\,\,=\,\, sinh^{-1}\Big(sin\,(u_i)\, sinh\,(t)\Big)
\tag{1}
\end{equation*}

\noindent for $i=1,2,3$. But for $t$ small we have $f_i(\alpha(t))=\frac{1}{2r}t_i(\alpha(t))$
and a simple differentiation of (1), evaluated at 0, shows $D(f_i)_p\, u =\frac{1}{2r}u_i$, $i=1,2,3$.
This proves (i), (ii) and (iii) and completes the proof of the proposition.\\

Choose $r>0$. In what follows we assume the width of the
normal neighborhoods of the $X_{\square}$ to be much larger than the number $r>0$. Lemma 2.1 asserts this is always possible.
Fix a point $p\in \dX_{\square^k}$ and consider $\Delta^{n-k-1}\sbs T_p X $. 
The cone $\rC_r \Delta^{n-k-1}$ is the set $\{tu\,,\, 0\leq t<r\, ,\, u\in\Delta^{n-k-1}\}$.
We have the exponential map 
$E:\rC_r\Delta^{n-k-1}\ra X$, given by $E(u,t)=exp_p(t\,u)$. 
Write $\square^n=\square^k\times\square^l\sbs\R^k\times\R^l=\R^n$ and
denote by $p\0{i}:\R^n\ra\R^i$ (with $i=k,l$), the projections onto the two factors.
Also, as in Section 6, we write $\bar{R}^n_+=[0,\infty)^n$.\\

\noindent {\bf Lemma 2.6.} {\it We have the following properties.}

\begin{enumerate}
\item[{(i)}] {\it The map $E$ respects faces, that is}\,\, {\small $E\,\Big(\,\, \big(\rC_r\Delta^{n-k-1}\big)
\,\,\,\cap \,\,\,T_p X_{\square^j}\,\,\Big)\,\, \,\sbs\,\,\, X_{\square^j}$}  
\item[{(ii)}]  {\it The map $f\circ E$ respects faces, that is}\,\, {\small $(f\circ E)\,\Big(\,\, \big(\rC_r\Delta^{n-k-1}\big)
\,\,\,\cap \,\,\,T_p X_{\square^j}\,\,\Big)\,\, \,\sbs\,\,\,\square^j$}. {\it Hence
$\big(\rC_r\Delta^{n-k-1}\big)
\,\,\,\cap \,\,\,T_p X_{\square^j}=\big(\rC_r\Delta^{n-k-1}\big)
\,\,\,\cap \,\,\,T_{f(p)}\square^j$.}  
\item[{(iii)}] {\it The map $p\0{l}\circ f\circ E$ does not depend on the point $p$.}
\item[{(iv)}]  {\it The values of the map $p\0{k}\circ f\circ E$ do not depend on the variable $u\in \Delta^{n-k-1}$ (but they do depend on $t$ and $p$).}
\item[{(v)}] {\it Write $T=(t_1,...,t_n)$.
The map $p\0{l}\circ T\circ E:\rC_s\Delta^{n-k-1}\ra\R^{n-k}_+$ is an embedding, provided $s<r-1$.}
\end{enumerate}

\noindent {\bf Remark.} In the second statement of (ii) we are considering both
$\rC_r\Delta^{n-k-1}\sbs T_pX $ and 
$\rC_r\Delta^{n-k-1}\sbs T_{f(p)}\square^n$.\\

\noindent {\bf Proof.} Statement (i) follows from the fact that each $X_{\square}$ is totally geodesic in $X$.
Statement (ii) also follows because $f$ respects faces. Differentiating,  using 
Proposition 2.5 and the fact that the derivative of $E$ (at 0) is the identity we obtain
$\big(\rC_r\Delta^{n-k-1}\big)
\,\,\,\cap \,\,\,T_p X_{\square^j}\sbs T_{f(p)}\square^j$. To get the other inclusion
use Proposition 2.5 again and count dimensions.\\

To prove (iii) and (iv) we assume for simplicity,   as in the proof of Proposition 2.5,  that
$k=n-3$ and $\square^{n-3}=\square_1^{n-1}\cap\square^{n-1}_2\cap\square^{n-1}_3$,
where $\square^{n-1}_i=\square^n\cap\{x_i=0\}$. Also $l=3$ and
$\square^3=\square_{4}^{n-1}\cap...\cap\square_n^{n-3}$. Let $u=(u\0{1}, u\0{2},u\0{3})\in\Delta^{2}$. %Hence $u_i=0$, $i>3$ and, 
Since $f=(f_1,...,f_n)$ we have 
$$ p\0{l}\circ f\circ E(u,t)\,\,=
\Big( f_1\big( E(u,t)  \big)      \, ,\,       f_2\big( E(u,t)  \big)      \, ,\, 
 f_3\big( E(u,t)  \big)      \,  \Big)
$$
\noindent But for $q\in X_{\square^{n-1}_i}\times [0,r)$ we have that
$f_i(q)=\frac{1}{2}\eta(t_i(q))$ (see paragraph before Proposition 2.4).
This together with (1) in the proof of 2.5 imply
$$  f_i\big( E(u,t)  \big)\, =\,\,\frac{1}{2}\,\eta\,\bigg(\,\, sinh^{-1}\Big(sin\,(u_i)\, sinh\,(t)\Big)  \,\,\bigg) 
$$
\noindent which depend only on $u$ and $t$. 
Here $u=(u_1,...,u_{n-k})$. This proves (iii).\\

Note that $p\0{k}\circ f\circ E(u,t)=\Big( f_4\big( E(u,t)  \big),..., f_n\big( E(u,t)  \Big)$.
We prove that  $f_4\big( E(u,t)  \big)$ is independent of $u$. The proof is similar for $i>4$. Write $U=X_{\square_4^{n-1}}$, $V=X_{\square^k}$, $W=U\cap V$ and let $p\in\dot{V}$. Note $U$, $V$, $W$ are
totally geodesic. Let $q\in W$ be such that the segment $[q,p]$ has length equal to the distance $d\0{X}(p, W)$ between
$p$ and $W$. 
If $d\0{X}(p, W)\geq r$ we are done because then $f_4\big( E(u,t)  \big)$ is constant. We assume
$L=d\0{X}(p, W)< r$.\\

Let $B$ be 
the union of the images of all geodesics of length $r$ in $U$ beginning at $q$ and perpendicular to $W$. Then
$B$ is isometric to the $r$-cone of the canonical all-right simplex $\Delta^{2}$, and it is totally geodesic.
Let $C$ be the union of the images of all geodesics of length $r$ in $X$ beginning at some point in $B$ and perpendicular to $U$. Then $C$ is isometric to $B\times [0,r)$ with the usual
$cosh$-warped product metric. We write $C=B\times [0,r)$.
Therefore $C$ is also totally geodesic, and $C$, $V$ intersect perpendicularly 
at $C\cap V=[q,p]$. Now, if $\alpha$ is a geodesic of length $<r$ beginning at $p$ and perpendicular to
$V$ then $\alpha$ in contained in $C$. Moreover, $\alpha$ is contained in $\ell\times [0,r)\sbs C$
for some ray $\ell\sbs B$ beginning at $q$. 
(To see this note that we can consider $C$ convex in
$\HH^n$, and the statement is true in $\HH^n$.)
Note that $\ell\times [0,r)\sbs C$ is isometric to a convex set in $\HH^2$. Finally
we get that $d\0{X}\big(\alpha(t),U\big)$ can then be computed in $\HH^2$ as the length of the side  $a$ of
a quadrilateral with consecutive  sides $a,b,c,d$, angles $\angle ab=\angle bc=\angle cd =\pi/2$ and $length(a)=t$, $length(b)=L$.
This calculation only depends on $t$ and $L=length [q,p]$ (hence on the choice of $p$) but not on the ``direction" $u$.\\

To prove (v) note that from equation (1) in the proof of Proposition 2.5 and the fact that $s<r-1$ we get that for $u=(u_1,...,u_{n-k})$, $0\leq t <s$, we have $t_i(tu)=sinh^{-1}\big(\,sin(u_i)\, sinh(t)\,\big)$.
This equation together with $\Sigma_{i=1}^{n-k}sin^2(u_i)=1$ imply 
$$t=sinh^{-1}\Big( \,\big(\,\Sigma_{i=1}^{n-k} \,sinh^2(t_i)\,\big)^{1/2}\,\Big).$$
\noindent Since we also get $u_i=sin^{-1}\big(\,  \frac{sinh(t_i)}{sinh(t)} \,\big)$,
the map $p\0{l}\circ T \circ E= (t_1\circ E,...,t_{n-k}\circ E)$ has a continuous inverse. Moreover, this inverse is clearly smooth when $t\neq 0$ and all $t_i<t$. But for $t=0$ the derivative of $p\0{l}\circ T\circ E$ can be shown to be injective.
This proves the lemma.\\

\noindent {\bf Remark 2.7.} Using the method in the proof above together
with hyperbolic trigonometry, we can find an explicit formula for the coordinate functions of the function in  (iv). It can be checked that
these maps are even in the variable $t$. 
In particular, $\frac{d}{dt}(p_k\circ f\circ E)|_{t=0}=0$.

\vspace{.5in}

\begin{center} {\bf \large 3. The Charney-Davis Hyperbolization Process.}
\end{center}

The strict hyperbolization process of Charney and Davis is done by gluing copies
of $X^n$ using the same pattern as the one used to obtain the cube complex $K$
from its cubes. This space is called $K_X$ in \cite{ChD}. 
We call this space the {\it piece-by-piece strict hyperbolization of $K$}.
Note that we get a map
$F:K_X\ra K$, which restricted to each copy of $X$ is just the map $f:X^n\ra\square^n$
in Section 2. We will write $X_{\square^k}=F^{-1}(\square^n)$, for a $k$-cube $\square^k$
of $K$.\\

But to obtain good differential and tangential properties, the process described above
is not enough. Therefore in \cite{DJ} and \cite{ChD} an alternative method is
given. We describe this next. As before 
let $X^n$ be a strict hyperbolization piece
and $K$ be a cube complex. We assume there is projection
$p:K\ra \square^n$ (see 7.2 of \cite{ChD}). Now consider $K_X$ given as the fiber product\\

$$\begin{array}{ccc}
K_X&  \stackrel{q\0{X}}{\longrightarrow}& X\\ \\
q\0{K}\downarrow&&\downarrow f\\ \\
K&  \stackrel{p}{\longrightarrow}& \square^n
\end{array} $$\\

\noindent that is $K_X=\{(y,x)\, :\, p(y)=f(x)\}\sbs K\times X$. Here $q\0{K}$, $q\0{X}$ are projections.
We call this space the {\it fiber-product strict hyperbolization
of $K$}.  We denoted both hyperbolizations by $K_X$ but we shall write $K_X^{{\mbox{\tiny piece-by-piece}}}$ and $K_X^{{\mbox{\tiny fiber-product}}}$
if we need to. We shall write $X_{\square^k}=q\0{K}^{-1}(\square^k)$, for a $k$-cube
$\square^k$ of $K$.\\

\noindent {\bf Remark.} The space $K_X^{{\mbox{\tiny fiber-product}}}$
does depend on the projection map $p$. For instance, if $p$ is a cube map,
that is $p|\0{\square^n}$ is an isometry for every $\square^n\in K$,  then $K_X^{{\mbox{\tiny piece-by-piece}}}$
and $K_X^{{\mbox{\tiny fiber-product}}}$ coincide. In general, these two
hyperbolizations are homeomorphic but the obvious homeomorphism (see below)
does not preserve the natural piecewise differentiable structures.
If needed we shall write $K_X^{{\mbox{\tiny fiber-product}}}(p)$
to show explicitly the dependence on $p$.\\

We now assume that $K$ has a smooth structure $\cS$ compatible with the
cube structure of $K$ (hence $1_K$ is a smooth cubulation of the smooth manifold $K$). We assume further that the 
projection $p:K\ra \square^n$ is smooth. Using this and item (4) at the beginning of Section 2 it is argued in \cite{ChD} that 0 is a regular value
of the smooth map $(k,x)\mapsto p(k)-f(x)$. Therefore $K_X$ is a smooth submanifold of $K\times X$
(with trivial normal bundle).
Hence if $K^n$ has a smooth structure (compatible with the cube structure $K$) then $K_X$ has a natural smooth structure.  
This is an important point for us, so we need to analyze this in a bit more detail.
First we remark two facts:

\begin{enumerate}
\item[{\bf (i)}] For the regular value argument to work it is assumed (implicitly) in \cite{DJ}, 1C.5, that
the restriction
$p|_{\dsquare}$ of $p$ to every open cube $\dsquare$ of $K$ is an embedding.
(In \cite{DJ} simplices are used instead of cubes.)
\item[{\bf (ii)}]  Let $\square$ be
a $k$-cube of $K$. Then $p(\square)=\square^k$, for some $k$-face $\square^k$ of $\square^n$.
If $K$ has no boundary the smoothness of $p$ implies that  for every $y\in \dsquare$ we have
 $Dp_y(T_yK)=T_{p(y)}\square^k$, where
$Dp$ is the derivative of $p$. 
In particular, the image  of $Dp_y$ is $k$-dimensional,
thus $Dp_y$ is not an isomorphism; hence $p|_\square$
is not an embedding. What is happening here is that
whenever we ``fold" two $n$-cubes into one (this is what $p$
does), and we want this folding $p$ to be smooth, then $p$
has to ``slow down to 0" at the place of the folding (which is
an $(n-1)$-cube). Of course none of this has to happen if, for
instance, the complex $K$ is equal to $\square^n$.)
\end{enumerate}

It is important for us here to work with both hyperbolization
constructions: $K_X^{{\mbox{\tiny piece-by-piece}}}$
and $K_X^{{\mbox{\tiny fiber-product}}}$. Therefore we
need a good way to identify them. We deal with this issue next. Note that $K_X^{{\mbox{\tiny piece-by-piece}}}$
has a natural piecewise smooth structure:
the inclusion of each copy of $X$ in
$K_X^{{\mbox{\tiny piece-by-piece}}}$ is,
by definition, a smooth embedding.  But this is not true
for $K_X^{{\mbox{\tiny fiber-product}}}$. We explain this next.\\

The {\it copy of $X$ in $K_X^{{\mbox{\tiny fiber-product}}}$  over an $n$-cube $\square^n$ of $K$} is $X_{\square^n}=q\0{K}^{-1}(\square^n)$. It is a ``copy" of $X$ because
the projection $q\0{X}|\0{X\0{\square^n}}:X_{\square^n}\ra X$ is a
homeomorphism, whose inverse is given by
$$
x\,\,\mapsto\,\, \Big(\, (p|\0{\square^n})^{-1}(\,  f(x)\,)\, ,\,x\,\Big).\,\,$$

\noindent But this map is not smooth because, as mentioned before, $p|_\square$ is not an embedding (we have
$Dp|_\square.v=0$, for some $v\neq 0$). Therefore, even though $q\0{X}$ is smooth, the map $q\0{X}|\0{X\0{\square^n}}$ is  not a diffeomorphism because
the natural (topological) embedding $\big(q\0{X}|\0{X\0{\square^n}}\big)^{-1}:X
\ra K_X$ is not smooth.\\

The price we paid for slowing down the cubes at the boundary
(via $p$)
was that we sped up the copies of $X$ at the boundaries.
Therefore the natural piecewise hyperbolic (and piecewise differentiable)
structure of $K_X^{{\mbox{\tiny piece-by-piece}}}$ does not directly
give one in $K_X^{{\mbox{\tiny fiber-product}}}$. \\

To have a chance to solve this problem we
need a more concrete expression for $p$. We will consider maps $p$ of the form $p=\bar{\rho}\circ c$, where $c$ is a 
{\it cube map} $c:K\ra \square^n$ (i.e. $c|\0{\square^n}$ is an isometry
for every $\square^n$) and $\bar{\rho}$ is
a {\it slow-down-at-the-boundary} map $\bar{\rho}:\square^n\ra\square^n$ given by $\bar{\rho}(x_1,...,x_n)\,=\, (\rho(x_1),...,\rho(x_n))$, with $\rho:I\ra I$ a smooth homeomorphism that
is a smooth diffeomorphism on $(0,1)$ and $\frac{d^k}{dt^k}\rho(0)=\frac{d^k}{dt^k}\rho(1)=0$, $k>0$. In what follows we write 
$K_X^{{\mbox{\tiny fiber-product}}}=K_X^{{\mbox{\tiny fiber-product}}}(\bar{\rho}\circ c)$.\\

\noindent {\bf Remark.} Let $c:K\ra\square^n$ be a cube map. For any $\square^n\in K$ the map
$(c|\0{\square^n})^{-1}:\square^n\ra\square^n\sbs K$ can be identified with
the inclusion $\square^n\hookrightarrow K$. In particular,
if $K$ has a smooth structure compatible with the
cube structure of $K$,  the map $(c|\0{\square^n})^{-1}$
is an embedding. Also, note that cube maps $c:K\ra\square^n$ are not smooth
(unless $K=\square^n$).\\

The following proposition says that $\bar{\rho}:\square^n\ra\square^n$ 
can be covered by a homeomorphism $X\ra X$.\\

\noindent {\bf Proposition 3.1.} {\it We can choose $\rho:I\ra I$
so that there is a smooth homeomorphism $P:X^n\ra X^n$ such that $f\circ P=\bar{\rho}\circ  f$,
i.e. the following diagram commutes}

$$
\begin{array}{ccc}
X&\stackrel{P}{\longrightarrow}&X\\
f\downarrow&&\downarrow f\\
\square^n&\stackrel{\bar{\rho}}{\longrightarrow}&\square^n
\end{array}
$$

\noindent {\it Moreover, we can choose $P$ so that
its restriction to every open face $\dX_{\square^k}$ is an embedding.}\\

\noindent {\bf Remark.} With a bit extra work we can get $P$ to be $B_n$-invariant,
but this fact will not be needed.\\

Note that from the construction of the map $\bar{\rho}$ we have that
$D\bar{\rho}|\0{q}.w=0$ for every $\square \in \square^n$, $q\in\square$ and $w$ perpendicular
to $\square$. We have the following addition to Lemma 3.1.\\

\noindent {\bf Addendum to Proposition 3.1.} {\it We can choose $P$ in
Proposition 3.1 so that for any $\square\in \square^n$ we have that
$DP|\0{p}.v=0$ for every $p\in\square$ and $v$ perpendicular
to $X_\square$.}\\

The proposition and its addendum are proved in Appendix A. \\

Now we get a new embedding $X\ra X_{\square^n}\sbs
K_X^{{\mbox{\tiny fiber-product}}}$ given by
$\big(q\0{X}|\0{X_{\square^n}}\big)^{-1}\circ P$. This is the ``correct" embedding,
as the next proposition shows.\\

\noindent {\bf Proposition 3.2.} {\it Suppose $p=\bar{\rho}\circ
c$ is smooth.
Then  $\big(q\0{X}|\0{X_{\square^n}}\big)^{-1}\circ P:
X\ra X_{\square^n}\sbs K_X^{{\mbox{\tiny fiber-product}}}$ is a smooth embedding. Moreover, the following
diagram commutes for every $n$-cube $\square^n$ of $K$.}\\

$$
\begin{array}{ccccc}
X&\xrightarrow{\,\,{\mbox{\tiny $\big(q\0{X}|\0{X_{\square^n}}\big)^{-1}\circ P$}}\,\,}&X_{\square^n}&
\xrightarrow{\,\,{\mbox{\tiny inclusion}}\,\,}& K_X\\\\
f\downarrow&& q\0{K}|\0{X_{\square^n}}\downarrow \,\,\,\,\,\,\,\,\,\,\,\,\,\,\,\,\,\,&&
\,\,\,\,\,\,\downarrow  q\0{K}\\\\
\square^n&\xrightarrow{\,\,\,\,\,\,\,\,\,\,\,\,\,\,\,\,1\0{\square^n}\,\,\,\,\,\,\,\,\,\,\,\,\,\,\,\,}&\square^n&
\xrightarrow{\,\,{\mbox{\tiny inclusion}}\,\,}&K
\end{array}
$$\vspace{.2in}

\noindent {\bf Proof.} 
From the definition of $K_X^{{\mbox{\tiny fiber-product}}}$ and Proposition 3.1 we get that the following diagram commutes

\begin{equation*}
\begin{array}{ccccc}
K_X&\xrightarrow{\,\,\,\,\,\,\,\,q\0{X}\,\,\,\,\,\,\,\,}&X&
\xleftarrow{\,\,\,\,\,\,\,\,P\,\,\,\,\,\,\,\,}& X\\\\
q\0{K}\downarrow&&f\downarrow &&
\,\,\,\,\,\,\downarrow f\\\\
K&\xrightarrow{\,\,\,\,\,\,\,\,p\,\,\,\,\,\,\,\,}&\square^n&
\xleftarrow{\,\,\,\,\,\,\,\,\bar{\rho}\,\,\,\,\,\,\,\,}&\square^n
\end{array}
\tag{1}
\end{equation*}\vspace{.2in}

\noindent The commutativity of (1) together with 
$p=\bar{\rho}\circ c$, and the fact that $(c|_\square)^{-1}$
is the inclusion $\square\hookrightarrow K$ imply that the left square of the diagram in the statement
of Proposition 3.2 commutes. The right square commutes by definition.\\

Write $g=\big(q\0{X}|\0{X_{\square^n}}\big)^{-1}\circ P$.
We have that the map $g:X\ra K_X\sbs K\times X$ is smooth if and only if the 
coordinate maps $q\0{K}\circ g$, $q\0{X}\circ g$ are smooth. 
First we have $q\0{X}\circ g=q\0{X}\circ \big(q\0{X}|\0{X_{\square^n}}\big)^{-1}\circ P=P$, which is smooth. \\

\noindent From diagram (1) and $p=\bar{\rho}\circ c$ we get 
$$q\0{K}\circ g=q\0{K}\circ 
\big(q\0{X}|\0{X_{\square^n}}\big)^{-1}\circ P=(p|\0{\square^n})^{-1}\circ \bar{\rho}\circ f=
\big(c|\0{\square^n}\big)^{-1}\circ \bar{\rho}^{-1}\circ\bar{\rho}\circ f=
\big(c|\0{\square^n}\big)^{-1}\circ f$$
\noindent Since $\big(c|\0{\square^n}\big)^{-1}$ is just the inclusion
$\square^n\ra K$ we get that $q\0{K}\circ g$ is smooth. 
It remains to prove that $g$ is a smooth embedding. 
Since $q\0{X}$ is smooth  and $P^{-1}$ is smooth on $\dX$
we get that $g^{-1}$ is smooth on $\dX_{\square^n}$. 
Note that the same argument shows that the restriction of $g$ to any
$\dX_{\square^k}$ is an embedding. Therefore if $u$ is  a non-zero
vector tangent to some $\dX_{\square^k}$ then $Dg.u$ is non-zero
and tangent to the corresponding $k$-face $g(X_{\square^k})$.
If $u$ is a non-zero vector normal to $\dX_{\square^k}$ then, by Lemma 2.5(i), $Df.u$ is non-zero and normal to $\square^k\sbs \square^n$ and certainly $D\big((c|\0{\square^n})^{-1} \circ f\big).u $ is also
non-zero and normal to $\square^k$. 
But from diagram (1)
we have that $q\0{K}\circ g=(c|\0{\square^n})^{-1}\circ f$,
hence $Dg.u$ is non-zero. Moreover,
$q\0{K}$ sends $g(X_{\square^k})$
to $\square^k$, therefore $Dg.u$ is not tangent to $g(X_{\square^k})$.
This proves that $Dg$ is injective on every point of $\dX_{\square^k}$.
This proves the proposition.\\

We can now use $\big(q\0{X}|\0{X_{\square^n}}\big)^{-1}\circ P$  on each
copy of $X$ in $K_X^{{\mbox{\tiny piece-by-piece}}}$ and get a map
$\Phi:K_X^{{\mbox{\tiny piece-by-piece}}}\ra K_X^{{\mbox{\tiny fiber-product}}}$
that is a smooth embedding on each copy of $X$. Hence we can
consider $K_X$ as $K_X^{{\mbox{\tiny piece-by-piece}}}$ with the pulled
back (by $\Phi$) differentiable structure, or $K_X^{{\mbox{\tiny fiber-product}}}$
with the pushed forward piecewise hyperbolic structure. \\

\noindent {\bf Corollary 3.3.} {\it The following diagram commutes.}\\

$$
\begin{array}{lcl}
K_X^{{\mbox{\tiny piece-by-piece}}}&\xrightarrow{\,\,\,\,\,\,\,\,\Phi\,\,\,\,\,\,\,\,}&K_X^{{\mbox{\tiny fiber-product}}}\\\\
F\downarrow&&\downarrow q\0{K}\\\\
K&\xrightarrow{\,\,\,\,\,\,\,\,1_K\,\,\,\,\,\,\,\,}&K
\end{array}
$$
\noindent {\it Moreover $\Phi$ is a smooth embedding on each copy of $X$ in $K_X^{{\mbox{\tiny piece-by-piece}}}$.} \vspace{.2in}

Here is an important caveat. We showed how to 
identify $K_X^{{\mbox{\tiny piece-by-piece}}}$
and $K_X^{{\mbox{\tiny fiber-product}}}$ in a good
way, so that we can benefit from the different
properties of both constructions. But a key piece
was missing: this identification
was done under the assumption that $p=\bar{\rho}\circ c$
is smooth (see statement of Proposition 3.2). We do not
know how to prove that $p$ is smooth, because the
smooth structure on $K$, though $PL$-compatible with
$K$, could be quite arbitrary. (If $p$ is not of the form
$p=\bar{\rho}\circ c$ the problem of finding $\Phi$ seems to be even
harder.) But in our case this does not matter because
we will work with normal smooth structures.
The next result shows that $p=\bar{\rho}\circ c$ is smooth 
%not just on $K=(K,\cS)$ but also on
on $(K,\cS')$, where $\cS'$ is a normal smooth structure on $K$ for $K$ (see Section 1).\\

\noindent {\bf Proposition 3.4.} {\it Let $\cS'$ be a normal smooth structure on $K$ for  $K$. Then $p :(K,\cS')\ra\square^n$ is smooth.}\\

The proof is presented in appendix B.\\

Note that we also have that the restriction $p|\0{\dsquare^i}$ on every open cube is an embedding, because the inclusions $\dsquare^i\ra (K,\cS')$ are also embeddings
(see remarks before Theorem 1.1.1).
Therefore the regular value argument in \cite{ChD} (see items (i) and (ii) at the beginning of Section 3)
goes through and we get the following result.\\

\noindent {\bf Corollary 3.5.} {\it We have that $K_X^{{\mbox{\tiny fiber-product}}}$ is a smooth submanifold of 
$(K,\cS')\times X$, with trivial normal bundle.}\\

We denote by $K_X'$ the submanifold $K_X^{{\mbox{\tiny fiber-product}}}\sbs (K,\cS')\times X$ with its
induced smooth differentiable structure.\\

The proof of Proposition 3.2 also works if we replace $K_X^{{\mbox{\tiny fiber-product}}}$ by $K_X'$, but with one change: we have to substitute
$X_{\square^n}$ by $\dX_{\square^n}$ (these are {\sf open} faces), that is, the map $\dX\ra\dX_{\square^n}\sbs K_X'$
is an embedding (the key point in the proof is that the inclusion 
$(c|\0{\square^n})^{-1}:\square^n\ra (K,\cS')$ is not an embedding, but its
restriction to $\dsquare^n$ is). It follows that
$\dX_\square$ is a submanifold of $K_X'$, for every $\square\in K$. Therefore we obtain the following corollary.\\

\noindent {\bf Corollary 3.6.} {\it  The map 
$$\Phi:K_X^{{\mbox{\tiny piece-by-piece}}}\xrightarrow{\,\,\,\,\,\,\,\,\,\,}K_X'$$ is a smooth embedding on each copy of $\dX$ in $K_X^{{\mbox{\tiny piece-by-piece}}}$.} \vspace{.2in}

\vspace{.6in}

\begin{center} {\bf \large 4. Normal neighborhoods on Charney-Davis Hyperbolizations.}
\end{center}

Theorem 1.1.1 gives a normal smooth atlas and a normal smooth structure for a given smooth
cubulation of a smooth manifold.
In this section we will construct a similar  atlas on the Charney-Davis strict hyperbolization of $K$.
 In what follows of this section we will use the notation $K_X$ for $K_X^{{\mbox{\tiny piece-by-piece}}}$.
Recall that we are denoting by $K_X'$ the submanifold $K_X\sbs (K,\cS')\times X$ with its
induced smooth structure. Recall we have 
a map $\Phi: K_X\ra K_X'$ (see Section 3).
The normal smooth structure $\cS'$ on $K$ has a normal atlas
$\cA=\cA\big(  \cL   \big)$, where
$\cL=\big\{  h_{\square^k}   \big\}\0{\square^k\in K}$ is a
smoothly compatible set of link smoothings for $K$.
We assume that the normal bundle of any face of the hyperbolization piece $X$
has width larger than $s\0{0}>0$.  Choose $r$, with $3r<s\0{0}$. By Lemma 2.1 the number $s\0{0}$ can be taken as large as we want.\\

By Lemma 2.5 we can use the derivative of
the map $F:K_X\ra K$ (in a piecewise fashion) to identify $\sL(X_{\square^k},K_X)$ with $\sL(\square^k,K)$, where in both cases we consider the ``direction" definition of
link, that is, the
link $\sL(X_{\square^k},K_X)$ (at $p\in \dX_{\square^k}$) is the set of normal vectors to $X_{\square^k}$ (at $p$) and  the
link $\sL(\square^k,K)$ (at $q\in\dsquare^k$) is the set of normal vectors to $\square^k$ (at $q$). Hence we write
$\sL(X_{\square^k},K_X)=\sL(\square^k,K)$; thus the set of links for $K$ coincides
with the set of links for $K_X$.\\

Let $X_{\square^k}\sbs X_{\square^{n}}$ be a $k$-face of $K_X$,
contained in the copy $X_{\square^n}$ of $X$ over $\square^n$. For a non-zero vector $u$ normal to $X_{\square^k}$ at $p\in X_{\square^k}$, and pointing inside 
$X_{\square^n}$, we have that $exp_p(tu)$ is defined and contained in
$X_{\square^n}$, for $0\leq t< s\0{0}/|u|$. Recall that $h_{\square^k}:\bS^{n-k-1}\ra\sL(\square^k,K)=\sL(X_{\square^k},K_X)$ is the smoothing of the
link corresponding to $\square^k$.
In the Introduction we defined the map

$$H\0{\square^k}\,\,\,:\,\,\D^{n-k}\times \dX_{\square^k}
\,\,\,\,\longrightarrow\,\,\,\, K_X $$
\noindent given by

$$H\0{\square^k}(\,t\,v\,\, ,\,\,p\,)\,\,\,=\,\,\, exp\0{p}\,\Big(\,\,2r\,t\,\,h_{\square^k}(v)\,\,\Big),$$\vspace{.2in}

\noindent where $v\in\bS^{n-k-1}$ and $t\in [0,1)$.
For $k=n$ we have that $H\0{\square^n}$ is the inclusion $\dX_{\square^n}\sbs K_X$
(or we can take this as a definition).
Note that $H\0{\square^k}$ is a topological embedding because we are assuming
the width of the normal neighborhood of $X_{\square}$ to be larger than $s\0{0}>2r$.
We called a chart of the form of $H\0{\square^k}$ (for some link smoothing
$h_{\square^k}$) a {\it normal chart for the $k$-face $X_{\square^k}$}.
A collection $\Big\{ H\0{\square^k} \Big\}\0{\square^k\in K}$ of normal charts
is a {\it normal atlas}, and if this atlas is smooth (or $C^k$) the induced
differentiable structure is called a {\it normal smooth (or $C^k$) structure.} The following is the Main Theorem in the Introduction.\\

\noindent {\bf Proposition 4.1.} {\it The normal  atlas $\Big\{ H\0{\square^k} \Big\}\0{\square^k\in K}$ on $K_X$
is smooth.}\\

\noindent {\bf Proof.} 
Since we are assuming $\cA=\cA\big(\{h\0{\square}\}\big)$ smooth we get from Proposition 1.3.4 that
the set  of smoothings $\{h\0{\square}\}$ is smoothly compatible, that is,
the maps in (1.3.2) (or 1.3.3) are smooth embeddings. For $\square^k\sbs\square^j\in K$
these maps have domains $\D^{n-j}\times \big(\dsquare^j\cap\sL(\square^k,K))\big)$ (the second factor is denoted
by $\dsigma\0{S}^j$ in Section 1.3) and target space $\bS^{n-k-1}$. We remark that
in this definition (and in Section 1.3) we use the ``geometric" definition
of link, while here in Section 4 we are using the ``direction" definition of link. But using  (piecewise Euclidean or piecewise hyperbolic) exponential maps in $K$ or $K_X$ we can identify these definitions.
Therefore we can identify $\square^j\cap\sL(\square^k,K)$ with
$\dX_{\square^j}\cap\sL(X_{\square^k},K)$ (the links here are geometric).
This together with the fact that $\sL(\square,K)=\sL(X_\square,K_X)$
imply that we can obtain maps  in the $K_X$ case similar to the maps
in (1.3.2), and these maps have the same domains and target spaces.
Moreover, they coincide modulo a slight smooth change (see remark below). Therefore
the $K_X$ versions of (1.3.2) are also ``smoothly compatible".
Now, the proof that
$\big\{ H\0{\square^k} \big\}$ is smooth is similar to the proof that
$\cA=\big\{ h^\bullet\0{\square^k} \big\}$ is smooth (assuming 
$\{h\0{\square}\}$ is smoothly compatible) given in the proof of Proposition 1.3.4. This proves Proposition 4.1.\\

\noindent {\bf Remark.} There is only one adjustment that has to be made in the
proof given in Proposition 1.3.4 to be applied to the case of Proposition 4.1.
In Section 1.2 (see (1.2.1)) we identified $\rC\sL(\square^k,K)$ as a subset of 
$\sL(\square^j,K)$, $\square^k\sbs\square^j$, using the radial projection
$\Re$ described in Remark 1.2.2. In the hyperbolic case (for Proposition 4.1),
hyperbolic radial projection give a similar identification  (call it $\Re\0{\HH}$). Moreover, 
since ray structures are preserved, these two projections coincide in
directions and just differ on the length. This length in the cube case
is given in Remark 1.2.2. Using hyperbolic trigonometry the analogous
formula (using the same setting as in Remark 1.2.2) is given by:
$tan^{-1}\Big(\frac{tanh\,(\,\,d\0{K}(v,p)\,\,)}{sinh\,(2r)}\Big)$, which is a 
smooth function. Therefore if  $\{h\0{\square}\}$ is smoothly compatible
using the identifications $\Re$, then  $\{h\0{\square}\}$ is also smoothly compatible
using the identifications $\Re\0{\HH}$.\\\\

We will denote by $\cS\0{K_X}=\cS\0{K_X}\Big( \big\{ h_{\square}\big\} \Big)$ the smooth structure on $K_X=
K_X^{{\mbox{\tiny piece-by-piece}}}$ induced by 
the smooth atlas $\cA\0{K_X}=\Big\{ H\0{\square^k} \Big\}\0{\square^k\in K}$.
Note that $\cA\0{K_X}$ depends uniquely on the smoothly compatible set of
link smoothings $\cL=\{h\0{\square}\}\0{\square\in K}$ for $K$ (hence for $K_X$),
and to express this dependence we will sometimes write 
 $\cA\0{K_X}=\cA\0{K_X}(\cL)$.\\

\noindent {\bf Proposition 4.2.} {\it The map $\Phi:\big(K_X,\cS\0{K_X}\big)\ra K_X'$ is a $C^1$-diffeomorphism.}\\

The proof is a bit technical and it is given in appendix C.\\

Hence the atlas  $\Big\{ H\0{\square^k}\Big\}$ is a normal $C^1$-atlas for the smooth 
manifold $K_X'$. The following is a more detailed version of the addendum to the Main Theorem given in the Introduction.\\

\noindent {\bf Proposition 4.3.} {\it The smooth manifolds $\big(K_X,\cS\0{K_X}\big)$ and $ K_X'$ are smoothly diffeomorphic.}\\

\noindent {\bf Proof.} Just approximate the $C^1$-diffeomorphism $\Phi$ by a smooth diffeomorphism.

\vspace{.6in}

\begin{center} {\bf \large 5. Normal structures for Hyperbolized Manifolds with Codimension Zero Singularities.}
\end{center}

In this section we treat the case of manifolds with a one point singularity.
The case of manifolds with many (isolated) point singularities
is similar. \\

We assume the setting and notation of Section 1.5.
Let $K_X$ be the Charney-Davis strict hyperbolization of $K$.
Denote also by $p$ the singularity of $K_X$.
Many of the definitions and results given in Sections 2, 3, 4 still hold (with minor changes)
in the case of manifolds with a one point singularity:\\

\begin{enumerate}
\item[{\bf (1)}] Given a set of link smoothings for $K$ (hence for $K_X$) 
we also get a set of charts $H_\square$ as in Section 4.
For the vertex $p$ we mean the cone map 
$H_p=\rC h_p:\rC N\ra \rC L\sbs K_X$. We will also denote the restriction
of $H_p$ to $\rC N-\{o\0{\rC N}\}$ by the same notation $H_p$.
As in item (3) of Section 1.5 here we are identifying $\rC N-\{o\0{\rC N}\}$ with $N\times (0,1]$ with the product smooth structure obtained from
\s{some} smooth structure ${\tilde{\cS}}_N$ on $N$.
As before $\{H_\square\}_{\square\in K}$ is a {\it normal atlas} for $K_X$
(or $K_X-\{p\}$). A normal atlas for $K-\{p\}$ induces a {\it normal smooth
structure on $K_X-\{p\}$.}
\item[{\bf (2)}]  Again we say that
the smooth atlas  $\{H_\square\}$ (or the induced smooth structure, or the set $\{h_\sigma\}$) 
is {\it correct with respect to $N$} if $\cS_N$ is diffeomorphic to ${\tilde{\cS}}_N$.
\item[{\bf (3)}] Let the set $\cL=\{h_\square\}_{\square\in K}$ 
induce a smooth structure on $K-\{p\}$, hence $\cL$ is smoothly compatible 
(see item (3) of Section 1.5). As in Proposition 4.1 we get that $\{H_\square\}_{\square\in K}$
is a smooth atlas on $K_X-\{p\}$ that induces a normal smooth
structure $\cS_{K_X}$ on $K_X-\{p\}$. Moreover, from Theorem 1.5.1 we get that $\cS_{K_X}$ is correct with respect to $\cS_N$ when $dim\, N\leq 4$ (always) or 
when $dim\, N>4$, provided $Wh(N)=0$. Note that in this case we can take
the domain $\rC N-\{o\0{\rC N}\}=N\times (0,1]$ of $H_p$ with smooth
product structure $\cS_N\times \cS_{(0,1]}$.
\item[{\bf (4)}] It can be verified that a version of Proposition 4.3 also
holds in this case: $(K_X-\{p\},\cS_{K_X})$ smoothly embeds in  $(K-\{p\},\cS')\times X$
with trivial normal bundle.
\end{enumerate}

\vspace{.6in}

\noindent {\bf \large  Appendix A. Proof of Proposition 3.1 and its Addendum.}\\

As always we write $I=[0,1]$.
Recall that the function $\bar{\rho}$ is defined as $\bar{\rho}(x_1,...,x_n)=
(\rho(x_1),...,\rho(x_n))$, where $\rho:I\ra I$ is as in Section 3. We will
assume the following extra condition on $\rho$:\\

\noindent {\bf (A.1.)} \hspace{1.5in} $\rho(x)=x$ for $\delta\leq x\leq 1-\delta$\\

\noindent for some small $\delta>0$. 
Let $\square^{n-1}$ be an 
$(n-1)$- face of $\square^n=\{(x_1,...,x_n)\, ,\, 0\leq x_i\leq 1\}$. For simplicity  write $\square^n=\square^{n-1}\times I$, and
consider the vector field on $\square^n$, depending on $\square^{n-1}$, given by $V\0{\square^{n-1}}(x)=e_n=(0,...,0,1)$. This vector field is perpendicular
to $\square^{n-1}$ and generates the collar $\eta\0{\square^{n-1}}:\square^{n-1}\times I\ra\square^n$,  of $\square^{n-1}$ in $\square^n$
(which for the decomposition $\square^n=\square^{n-1}\times I$ is just the identity). \\
%Note that $V_{n-1}$ is invariant by the
%action of $B_{n-1}$ on $\square^n$, where $B_{n-1}=B_{\square^{n-1}}$ is the %subgroup of elements $b$ in $B_n$ such that $b(\square^{n-1})=\square^{n-1}$. \\

Let $\hat{\rho}$ be the smooth self-homeomorphism on $\square^{n-1}\times [0,\delta]$
given by $\hat{\rho}(x,t)=(x,\rho(t))$. 
Let $\Lambda\0{\square^{n-1}}$ be the smooth self-homeomorphism on
$\square^n\ra\square^n$ that is the identity outside $\eta\0{\square^{n-1}}\big( \square^{n-1}\times [0,\delta) \big)$  and on 
the image of $\eta\0{\square^{n-1}}$ it is equal
to $\eta\0{\square^{n-1}}\circ\hat{\rho}
\circ\eta\0{\square^{n-1}}^{-1}$. Hence we can write\\

\noindent {\bf (A.2.)} \hspace{1.8in} 
$\bar{\rho}\,\,=\,\, \Lambda\0{\square_1^{n-1}}\circ...\circ \Lambda\0{\square^{n-1}_{2n}}$\\

\noindent for any ordering $\square_1^{n-1},...,\square^{n-1}_{2n}$ of all the
$(n-1)$-faces of $\square^n$.\\

We will assume that the width of the normal neighborhoods of the
$X_{\square}$ in $X$ are larger than $3r$ (see Section 2).
\\

\noindent {\bf Lemma A.3.} {\it For each $\square^{n-1}$ the vector field
$V\0{\square^{n-1}}$ has a lifting $W_{\square^{n-1}}$ to $X$ near $X_{\square^{n-1}}$.
Moreover, $W_{\square^{n-1}}$
%is invariant by the action of $B_{n-1} on $X$ and it 
is perpendicular to $X_{\square^{n-1}}$.}\\

\noindent {\bf Remark.}
%\noindent {\bf 1.} 
By $W_{\square^{n-1}}$ being a lifting of $V_{\square^{n-1}}$
 near $X_{\square^{n-1}}$ we mean that $W_{\square^{n-1}}$ is
defined on a normal neighborhood of $X_{\square^{n-1}}$ of width $\leq r$,
and $Df\,.\, W_{\square^{n-1}}=V_{\square^{n-1}}$.\\

%\noindent {\bf 2.} Here $B_{n-1}=B_{\square^{n-1}}$ is formed by elements $b$ in %$B_n$ such that $b(X_{\square^{n-1}})=X_{\square^{n-1}}$. \\

Before we present the proof of Lemma A.3 we show how it implies
Proposition 3.1. The addendum to Proposition 3.1 will be proved later, at the end of
this appendix. There is an $s'$ such that all $W_{\square^{n-1}}$
are defined on the normal neighborhood of $X_{\square^{n-1}}$ of width $s'$.
Using the vector fields $W_{\square^{n-1}}$ we get collars
$\tau\0{\square^{n-1}}:X_{\square^{n-1}}\times[0,a]\ra X$, for some fixed $a>0$.
Since $W_{\square^{n-1}}$ is a lifting of $V_{\square^{n-1}}$
we get\\

\noindent {\bf (A.4.)}  \hspace{1.4in} 
$f\,\Big( \tau\0{\square^{n-1}}\big( x,t  \big) \Big)\,\,=\,\,
\eta\0{\square^{n-1}}\big( f(x),t  \big).$\\

\noindent For instance, in the special case of the trivial decomposition
$\square^n=\square^{n-1}\times I$ we get 
$f\,\big( \tau\0{\square^{n-1}}\big( x,t  \big) \big)=
\big( f(x),t  \big)$ because, in this case $\eta\0{\square^{n-1}}$
is just the identity.
Let  now $\theta\0{\square^{n-1}}$ be the smooth self-homeomorphism on
$X_{\square^{n-1}}\times[0,a]$ given by\\

\noindent {\bf (A.5.)}  \hspace{1.7in} 
$\theta\0{\square^{n-1}}\big( x,t  \big) \,\,=\,\,\big( x,\rho(t)  \big).$\\

Assuming $\delta>0$ in (A.1) such that $\delta<a$, we get that  $\theta\0{\square^{n-1}} $ is the identity outside 
$X_{\square^{n-1}}\times[0,\delta]\sbs X_{\square^{n-1}}\times[0,a)$.
Finally define $\Theta\0{\square^{n-1}}$ to be the the smooth self-homeomorphism on $X$ that is the identity outside $\tau\0{\square^{n-1}}\big( X_{\square^{n-1}}\times [0,\delta) \big)$  and on 
the image of $\tau\0{\square^{n-1}}$ is equal
to $\tau\0{\square^{n-1}}\circ\theta\0{\square^{n-1}}
\circ\tau\0{\square^{n-1}}^{-1}$.\\

\noindent {\bf Claim A.6.} {\it For every \,$\square^{n-1}$ we have that}\,\,\,
$ f\circ \Theta_{\square^{n-1}}\,\,=\,\,\Lambda_{\square^{n-1}}\circ f.$\\

\noindent {\bf Proof of Claim.} By (A.4) we have that
$f\big( \tau\0{\square^{n-1}}\big(X_{\square^{n-1}}\times [0,\delta)\big) \big)=
\eta\0{\square^{n-1}}\big(\square^{n-1}\times [0,\delta)\big)$.
Hence a point $p\in X$ is in $\tau\0{\square^{n-1}}\big(X_{\square^{n-1}}\times [0,\delta)\big)$ if and only if its image $f(p)$ is in
$\eta\0{\square^{n-1}}\big(\square^{n-1}\times [0,\delta)\big)$.
If $p$ is not in $\tau\0{\square^{n-1}}\big(X_{\square^{n-1}}\times [0,\delta)\big)$ we get that $\Theta_{\square^{n-1}}(p)=p$ and $\Lambda_{\square^{n-1}}(f(p))=f(p)$ and the claim is true in this case. Assume now that $p$ is in
$\tau\0{\square^{n-1}}\big(X_{\square^{n-1}}\times [0,\delta)\big)$.
Write $\tau\0{\square^{n-1}}(x,t)=p$, thus $f(p)=\eta\0{\square^{n-1}}( f(x),t)$.
By applying (A.4) and (A.5) several times we get

$$\begin{array}{ccl}
f\circ \Theta\0{\square^{n-1}}(p)&=&f\circ \tau\0{\square^{n-1}}\circ\theta\0{\square^{n-1}}\circ\tau\0{\square^{n-1}}^{-1}(p)\\
&=&
f\circ \tau\0{\square^{n-1}}\circ\theta\0{\square^{n-1}}\big(x,t\big)\\ &=&
f\circ \tau\0{\square^{n-1}}\big(x,\rho(t)\big)\\ &=& \eta\0{\square^{n-1}}\big(f(x),\rho(t)\big)
\\ &=& \eta\0{\square^{n-1}}\circ\hat{\rho}\big(f(x),t\big)
\\ &=& \eta\0{\square^{n-1}}\circ\hat{\rho}\circ\eta\0{\square^{n-1}}^{-1}\circ f(p)\\
&=&\Lambda\0{\square^{n-1}}\circ f (p).
\end{array}
$$

\noindent This proves the claim.\\

To finish the proof of Proposition 3.1 just define $P=\Theta_{\square^{n-1}_1}\circ...\circ\Theta_{\square^{n-1}_{2n}}$.
The fact that $f\circ P=\bar{\rho}\circ f$ follows from (A.2) and Claim A.6. This proves
Proposition 3.1. \\

It remains to prove Lemma A.3. Note that to prove Lemma A.3 it is not enough to use Lemma 2.2 and
a local diffeomorphism argument to obtain the lifting, because $f$ has singularities on the boundary
(for $n\geq 3$).\\

\noindent {\bf Proof of Lemma A.3.}
Fix $\square^{n-1}$. Without loss of generality we assume $\square^{n-1}=
\square^{n-1}_1$, where $\square^{n-1}_i=\{\,x_i=0\,\}\cap\square^n$.
We have $\square^n=I\times\square^{n-1}_1$. Write
$V=V_{\square^{n-1}_1}$ and $W=W_{\square^{n-1}_1}$. 
Now, since the condition $Df.W=V$ is linear we have, using a partition of
unity
%a normalization 
and taking $\delta$ small in (A.1), that it is enough to find
a lift of $V$ just locally, that is:\\

\begin{enumerate}\item[{\bf (A.7.)}] {\it  for every $p\in\square^{n-1}_1$
there is a neighborhood $U$ of $p$ in $X$ and vector field
$W$ on $U$ such that $Df.W=V$}.\end{enumerate}
%\vspace{.2in}

%\noindent (To obtain the invariance by $B_{n-1}$ condition take $\frac{1}%{N}\Sigma_{b\in B_{n-1}} Db.W$ for any $W$ with $Df.W=V$. Here %$N$ is the cardinality of $B_{n-1}$.)
%Before we prove (A.7) we need to introduce some objects.
%\\

Let $\square^k\sbs\square^{n}$. Write
$D^j(s)=\rC_s\Delta^j=\{tu\in\R^{j+1},\,t\in[0,s],
\,u\in\Delta^j \}$.
We identify the closed normal neighborhood 
$N_s(\square^k)$ of $\square^k$ 
of width $s$
with $\square^k\times D^{n-k}(s)$ (here $s<1$).
Similarly we identify the closed normal neighborhood 
$N_s(X_{\square^k})$ of $X_{\square^k}$ 
of width $s$
(via the exponential map) with $X_{\square^k}\times D^{n-k}(s)$.
Note that for $\square^k\sbs\square^n$ 
%let $A(\square^k)=\{\square^{n-1}\,,\,\square^k\sbs\square^{n-1}\}$. 
we can write $\square^k=\bigcap\0{\square^k\sbs\square^{n-1}}\square^{n-1}$.
Define

$$A_s(\square^k)=\bigcap\0{\square^k\sbs\square^{n-1}}N_{s}\Big(\square^{n-1}  \Big),$$

 $$A_s(X_{\square^k})=\bigcap\0{\square^k\sbs\square^{n-1}}N_{s}\Big(X_{\square^{n-1} } \Big)$$
\noindent and  for $k<n$

$$
\begin{array}{ccl}L(X_{\square^k})&=& A_{3r}\big(\, X_{\square^k}   \,  \big)
\,\,\,-\,\,\, \bigcup\0{\square^k\not\subset\square^{n-1}}N_{2r}\Big(X_{\square^{n-1}}  \Big)\\
\\
&=&
\bigcap\0{\square^k\sbs\square^{n-1}}N_{3r}\Big(X_{\square^{n-1}}  \Big)\,\,\,-\,\,\, \bigcup\0{\square^k\not\subset\square^{n-1}}N_{2r}\Big(X_{\square^{n-1}}  \Big).
\end{array}$$\\

Note that $A_{s}(X_{\square^k})\sbs N_{s'}(X_{\square^k})$ for large $s'$
(how large $s'$ should be with respect to $s$ can be calculated using hyperbolic
trigonometry). Hence $L(X_{\square^k})\sbs N_{s}(X_{\square^k})$ for large $s$.\\

%and $L(X_{\square^n})=X_{\square^n}-\bigcup\0{\square^{n-1}\sbs\square^n}N_{2r}
%\Big(X_{\square^{n-1}}  \Big)$.\\

\noindent {\bf Claim A.8.} {\it We have 
 $X_{\square^{n-1}_1}\sbs\bigcup\0{\square^k\sbs \square^{n-1}_1}
L(X_{\square^k})$.}\\

\noindent 
Let $p\in X_{\square^{n-1}_1}$. 
If $p\notin L(X_{\square^{n-1}_1})$ then $p\in N_{2r}(X_{\square^{n-1}})$ for some
$\square^{n-1}$. Hence $p\in A_{2r}\big( X_{\square^{n-1}\cap\,\square^{n-1}_1}\big)$. Therefore we either have $p\in   L_{2r}\big( X_{\square^{n-1}\cap\square^{n-1}_1}\big)$, or $p\in N_{2r}(X_{\square_2^{n-1}})$,
for some $\square^{n-1}_2$ different from $\square^{n-1}_1$ and $\square^{n-1}$.
Arguing in the same way by induction we get that if $p\notin L(X_{\square^k})$
for all $\square^k\sbs\square^{n-1}_1$ with $k>0$ then $p\in A_{2r}(X_{\square^0})\sbs L(X_{\square^0})$, for some vertex $\square^0$. This proves Claim A.8.\\

We now prove statement (A.7).
We use the construction of the map $f$ given in Section 2.
Let $p\in \square^{n-1}_1$.
From Claim A.8 we can assume that $p\in L(X_{\square^k})$, for some $\square^k\sbs\square^{n-1}_1$. Write $l=n-k$.
Note that  $L(X_{\square^k})\sbs N_{s}(X_{\square^k})=X_{\square^k}
\times D^l(s)$ (for large $s$), hence we will sometimes write $p=(p,0)\in X_{\square^k}
\times D^l(s)=N_{s}(X_{\square^k})$.\\

For simplicity we assume $\square^k=\square^{n-1}_{1}\cap...\cap\square^{n-1}_{l}$, $l=n-k$. Hence, using the notation in Lemma 2.6, we have that 
$p\0{l}\circ f=(f_1,...,f_l)$ and $p\0{l}\circ T=(t_1,...,t_l)$.\\

\noindent {\bf Claim A.9.} 
{\it We have }
\begin{enumerate}
\item[{\it (a)}] {\it if $i >l$ then $f_i(p)=1/2$,}
\item[{\it (b)}] {\it if $(q,u)\in L(X_{\square^k})$ is
close to $p=(p,0)$, then $f_i(q,u)=1/2$, $i>l$,}
\item[{\it (c)}] {\it let $U=U'\times D\sbs X_{\square^k}\times D^l(s)$
be a product neighborhood where (b) holds for every $(q,u)\in U$.
Then $p\0{l}\circ T $ is an embedding on $ \{q\}\times D$, for every $q\in U'$.}
\end{enumerate}

\noindent Since $p\in L(\square^k)$, we have that $p\notin N_{2r}(X_{\square^{n-1}_i})$,
for $i>l$. Therefore $t_i(p)=d\0{X}(p,X_{\square^{n-1}_i})>2r>r$, $i>l$, and (a)
follows. Item (b) follows from (a), continuity and the fact that
the sets $N_s$ are closed. Item (c) follows from
Lemma 2.6 (v). This proves Claim A.9.\\

To finish the proof of (A.7) on $U$ just take $W(q,u)=\frac{1}{2r}\Big(\big( p\0{l}\circ T \big)|\0{\{      q\}\times D}\Big)^*(e_1)$, where $e_1$ is the constant vector field
$(1,0,...,0)$ on $\R^l$. (Note that $W$ is different from the gradient, with respect to
the hyperbolic metric on $X$, of the distance to $X_{\square^{n-1}_1}$ function $t_1$.)
It follows now from (b) of Claim A.9 and the fact that $f_1(x)=\frac{1}{2r}\rho(t_1(x))=\frac{t_1(x)}{2r}$, if $x$ is close to $X_{\square^{n-1}_1}$, 
that $Df.W=e_1=V$. This proves (A.7). It can be verified from the
construction that the second statement of A.3 holds.  This proves Lemma A.3 and completes the proof of Proposition 3.1. \\

\noindent {\bf Proof of the Addendum to Proposition 3.1.}
%For each $\square^{n-1}$ let {\bf n}$\0{\square^{n-1}}(p)$ be the normal
%vector to $X\0{\square^{n-1}}$ at $p\in X\0{\square^{n-1}}$.
%From A.3 and 2.5 (ii) we have $V\0{\square^{n-1}}=2r{\mbox{\bf n}}$.
Since for $p\in\square^{n-1}$ we have that $W_{\square^{n-1}}(p)$ is perpendicular to $X_{\square^{n-1}}$,
it is enough to prove that $DP_p.W_{\square^{n-1}}(p)=0$,
for every $\square^{n-1}$ and $p\in\square^{n-1}$.
To make this happen we need to modify our construction
of $P$ a little bit.\\

Note that from the second statement in Lemma A.3, (A.5) and the definition of
$\Theta\0{\square^{n-1}}$ we get\\

\noindent {\bf (A.10.)}\hspace{1.5in}
$D\Theta\0{\square^{n-1}}.W\0{\square^{n-1}}\,=\,0$\\

We need a lemma, which is essentially an initial value version of Lemma A.3.\\

\noindent {\bf Lemma A.11.} {\it Let $U$ be a (not necessarily tangent) vector field on $X\0{\square^{n-1}}$. Suppose that
$Df.U=V\0{\square^{n-1}}$. Then there is a self-diffeomorphism
$g$ on $X$ covering the identity $1\0{\square^{n}}:\square^{n}\ra\square^{n}$
(see diagram) with $Dg.U=W\0{\square^{n-1}}$.}
$$
\begin{array}{ccc}
X&\stackrel{g}{\longrightarrow}&X\\
f\downarrow&&\downarrow f\\
\square^{n}&\stackrel{1\0{\square^n}}{\longrightarrow}&\square^n
\end{array}
$$

\noindent {\bf Proof.} Using collars and integral curves the problem is reduced to
finding an extension $U$ of $U$ to a neighborhood of $X\0{\square^{n-1}}$,
with $Df.U=V\0{\square^{n-1}}$ and $U=W\0{\square^{n-1}}$
outside an even smaller neighborhood of $X_{\square^{n-1}}$
(the argument uses the integral curves of $-U$).
The proof that such an extension exists is similar to that of Lemma A.3.
(without the perpendicularity condition). The only change needed is at the very end of the proof of (A.7) (after the proof of (A.9)). In our present case 
we have that $U(q)=U(q,0)=W(q,0)+T(q)$, where $T(q)$ is tangent
to $X\0{\square^{n-1}}$ (this is because $Df.U=V\0{\square^{n-1}}$ and $q\in L(X\0{\square^{k}})$).
Now take
$U(q,v)=W(q,u)+\rho(|v|)T(q)$, where $\rho(t)$ is equal to 1 near $t=0$ and equal to
0 for $t\geq \mu$, for some small $\mu>0$. This proves the lemma.\\

We now prove the addendum. 
Recall that at the beginning of Appendix A we ordered the $(n-1)$-cubes: $\square^{n-1}_1,...
,\square^{n-1}_{2n}$, and we constructed the corresponding
$\Lambda_{\square^{n-1}_i}$, $\Theta_i=\Theta_{\square^{n-1}_i}$.
Write $V_i=V\0{\square^{n-1}_i}$
and  $W_i=W\0{\square^{n-1}_i}$.
We will need the following statement
which follows from the definition of the $\Lambda_{\square^{n-1}_i}$.\\

\noindent {\bf (A.12.)}\hspace{1.5in} $D\Lambda_{\square^{n-1}_i}.V_j=V_j$,
{\it \,\,\,for $i\neq j$}.\\

Now take $P=\Theta_{2n}\circ g\0{2n-1}\circ...\circ g\0{1}\circ\Theta_1$, where
the $g\0{i}$ are obtained in the following way. 
From Claim A.6, (A.12) and the fact that $Df.W_i=V_i$ we get that $Df(D\Theta_1. W_2)=V_2$, hence we can apply Lemma
A.11 to get a self-diffeomorphism $g\0{1}:X\ra X$ lifting the identity and 
satisfying $Dg\0{1}.(D\Theta_1.W_2)=W_2$. Next note that from
 A.6, (A.12), Lemma A.11 and the fact that $Df.W_i=V_i$ we get that 
$Df(D(g\0{1}\circ \Theta_1). W_3)=V_3$ and 
we can apply Lemma
A.11 to get a self-diffeomorphism $g\0{2}:X\ra X$ lifting the identity and 
satisfying $Dg\0{2}.(D(g\0{1}\circ\Theta_1).W_3)=W_3$, and so on.
From the choice of the $g\0{i}$ and (A.10) we get that $DP. W_i=0$.
Also from Claim A.6 and the fact that $g\0{i}$ lifts the identity we
get $f\circ P=\bar{\rho}\circ f$. This proves the addendum to Proposition 3.1.
\vspace{.6in}

\noindent {\bf \large  Appendix B. Proof of Proposition 3.4.}\\

We shall demand the following condition on $\rho$: that the derivatives of $\rho$ {\it approach zero exponentially fast, at 0 and 1}. That is\\

\begin{enumerate}
\item[{\bf (B.1.)}]\hspace{.3in} {\it for every $k$ there are positive $a$ and $b$ such that $|\frac{d^k}{dt^k}\rho(t)|\leq a e^{-\frac{b}{(1-t)t}}$.}
\end{enumerate}\vspace{.2in}

  Let $\cA=\Big\{ \big(   h_{\square^i}^\bullet,\D^{n-i}
\times\dsquare^i    \big)  \Big\}$ be a normal atlas inducing $\cS'$.
We write $W_{\square^i}$ for the image of $h^\bullet_{\square^i}$. Note that
$W_{\square^i}$ is a normal neighborhood $\rC\sL (\square, K)\times\dsquare^i$
of $\dsquare^i$.
Write $c=(c_1,...,c_n)$.
We will prove that $\mu=\rho\circ c_1$ is smooth. The proof for $\rho\circ  c_i$ is the same.\\

We state three facts about the map $c_1$, which can be verified by inspecting
each of them cube by cube.

\begin{enumerate} 
\item[{\bf (1)}] There are three possibilities for a cube $\square\in K$:
First  $c_1(\square)=\{0\}$ and we say $\square$ is a {\it 0-valued-cube},
second $c_1(\square)=\{1\}$ and we say $\square$ is a {\it 1-valued-cube}
and finally $c_1|_{\square}$ is onto $I=[0,1]$ and we say in this last
case that $\square$ is an {\it I-cube}.
In what follows everything we do for 0-valued-cubes can be done for 1-valued-cubes, so we will
just ignore 1-valued-cubes.
\item[{\bf (2)}] For a 0-valued-cube $\square^i$ the map $c_1$ (and hence $\mu$ and all
its derivatives)
is a product map on a neighborhood of $\square^i$. Specifically $c_1$ factors
through a composition
$$
W_{\square^i}\,\,=\,\, \rC\sL (\dsquare^i,K)\times \square^i\,\,
\stackrel{{\mbox{{\tiny projection}}}}{\longrightarrow} \,\,\rC\sL (\dsquare^i,K)
\,\,\longrightarrow \,\, I.
$$
\item[{\bf (3)}]
For a I-cube $\square^i$ the map $c_1$ (and hence $\mu$ and all its derivatives)
is a product map on a neighborhood of $\square^i$. Specifically $c_1$ factors
through a composition
$$
W_{\square^i}\,\,=\,\, \rC\sL (\dsquare^i,K)\times \dsquare^i\,\,
\stackrel{{\mbox{{\tiny projection}}}}{\longrightarrow} \,\,\square^i
\,\,\longrightarrow \,\, I,
$$
\noindent where the last arrow is also a projection: $(x_1,...,x_n)\mapsto x_1$.
\end{enumerate}\vspace{.2in}

We prove that $\mu=\rho\circ c_1$ is smooth by showing that its representative
$\mu_{\square}=\mu\circ (h_{\square}^\bullet)^{-1}$ on each chart is smooth. We prove this by induction
on the decreasing dimension of the cubes. Consider first the following two statements
that depend on the $i$-cube $\square^i$:\\

\begin{enumerate}
\item[{\bf A($\square^i$):}]  We have that $\mu_{\square^i}$ is smooth on $\D^{n-i}\times\dsquare^i$. (Hence $\mu$ is smooth on $W_{\square^i}$.)
\item[{\bf B($\square^i$):}]  For every $\square^j<\square^i$, $\square^j$ a 0-valued-cube,
the map $\mu_{\square^i}$ and all its derivatives approach zero
exponentially fast with respect to the distance to $\square^j$. That is,
 for every $k$ there are positive $a$ and $b$ such that $|\frac{d^k}{dt^k}
\mu_{\square^i}(p)|\leq a e^{-\frac{b}{(1-t)t}}$, where $t=d\0{\D^{n-i}\times\square^i}(p,\square^j)$, $p\in \D^{n-i}\times\dsquare^i$.
\end{enumerate}

Recall that the chart maps $h_{\square^i}^\bullet$ respect the product $\D^{n-i}
\times\dsquare^i$ and the inclusion maps of cubes $\dsquare\ra (M,\cS')$ are
embeddings.
% (hence $\dsquare\ra (M,\cS')\stackrel{q}{\ra}\square^n$ is the inclusion). 
Therefore item (3) above implies that for an $I$-cube $\square^i$ the map 
$\mu_{\square^i}$ is just projection given by the composition
$$
\D^{n-i}\times\dsquare^i\stackrel{{\mbox{{\tiny projection}}}}{\longrightarrow}
\dsquare^i\stackrel{{\mbox{{\tiny projection}}}}{\longrightarrow}\,\,\, I\,\,\,
\stackrel{\rho}{\longrightarrow}\,\,\, I,
$$
\noindent where the last projection is projection to the $x_1$ coordinate.
Since $\mu_{\square^i}=\rho\circ \pi$, where $\pi$ is linear,
we have that {\bf A($\square$)} is true for every $I$-cube $\square$.
Also, if $\square^i$ is an $I$-cube and $\square^j<\square^i$ is a 0-valued-cube,
we can write $\square^i=\square^{i-1}\times \square^1$, $\square^j<\square^{i-1}$,
where the projection on to the $x_1$ coordinate is $\square^{i-1}\times \square^1
\ra \square^1=I$. But for $p\in \D^{n-i}\times\dsquare^i$ we have
$$\pi(p)\,\,\,\leq\,\,\,d\0{\D^{n-i}\times\square^i}(p,\square^{i-1})\,\,\,\leq\,\,\, d\0{\D^{n-i}\times\square^i}(p,\square^j).$$
\noindent This together with $\mu_{\square^i}=\rho\circ\pi$,
(B.1), and the fact that $\pi$ is linear imply that {\bf B($\square^i$)} is also true for every $I$-cube $\square^i$.\\

For a 0-valued $(n-1)$-cube $\square^{n-1}$ it is straightforward to verify that
{\bf A($\square^{n-1}$)} and {\bf B($\square^{n-1}$)} hold true.
Assume now that {\bf A($\square^i$)} and {\bf B($\square^i$)} hold true for every
0-valued-cube $\square^i$, $i> k$. We prove the same is true for 0-valued-cubes $\square^k$.\\

Let $\square^k$ be a cube of dimension $k$.
Since $W_{\square^k}-\dsquare^k\sbs \bigcup_{i>k}W_{\square^i}$
it follows from the inductive hypothesis
{\bf A($\square^i$)}, $i>k$, that $\mu_{\square^k}$ is smooth on $\D^{n-k}\times\dsquare^k- \dsquare^k$, where we are writing $\dsquare^k=\{0\}\times \dsquare^k$.
By item (2) above $\mu_{\square^k}$ is a product on $\D^{n-k}\times\dsquare^k$,
hence,  it is enough to prove that the restriction $\nu=\mu\0{\square^k}|_{\D^{n-k}}:\D^{n-k}\ra I$
is smooth at  $0\in\D^{n-i}$. And
by the Mean Value Theorem we only need to prove that all partial
derivatives of $\nu:\D^{n-k}\ra I$ tend to zero as
a point $p$  tends to $0\in\D^{n-i}$.\\  

Let $v_m=t_mu_m\in\D^{n-k}$,  $u\in\bS^{n-k-1}=\p \D^{n-k}$, $t_m\in(0,1)$, $t_m\ra 0$.
We want to prove that all partial derivatives of $\nu$ at $v_m$ tend to zero
as $m\ra \infty$.
We can assume (arguing by contradiction) that  $u_m\ra u\in\bS^{n-k-1}$.\\

Corollary 1.1.2 says that the link $S=\sL(\square^k,K)$ is
a submanifold of $(M,\cS')$. The open sets 
$U_{\square^i}=S\cap W_{\square^i}$, $\square^i>\square^k$, form an open cover of $S$. 
Note that $U_{\square^i}$ is a normal neighborhood of $\square^i\cap S$
in $S$.
Write $x_m=h_{\square^k}^\bullet (v_m)$ and $y_m=h_{\square^k}^\bullet (u_m)\in S$. (Rigorously $h_{\square^k}^\bullet$ is not defined on $\p\D^{n-k}$ but,
after rescaling, we can assume this does happen).
Since $h_{\square^k}^\bullet$ restricted to $\D^{n-k}$ is, by definition, a cone map 
we can write $x_m=t_my_m$, where this last product is realized on the cone
link $\rC\sL (\square^k,K)$ of $\square^k$. And we also get $y_m\ra z=h_{\square^k}^\bullet
(u)$. \\

We
% can assume (again arguing by contradiction) 
have that $z\in\dsquare^i\cap S$, for some $\square^i>\square^k$. Let $V$ be a small neighborhood of $z$ in $S$ with
$\bar{V}\sbs U_{\square^i}$ and we assume $y_m\in V$ for all $m$.
Write $$\nu=\mu\circ \big(h_{\square^k}^\bullet|_{\D^{n-i}}\big)= \Big(\mu\circ h_{\square^i}^\bullet\Big)\circ\Big( \big(h_{\square^i}^\bullet\big)^{-1}
\big(h_{\square^k}^\bullet|_{\D^{n-i}}\big)\Big).$$

By {\bf B($\square^i$)} all partial derivatives of the first term
$\mu\circ h_{\square^i}^\bullet$ approach zero exponentially fast as a point
get close to $\square^k$. Likewise, by Corollary 1.4.3 the derivatives of the second term
$ \big(h_{\square^i}^\bullet\big)^{-1}
\big(h_{\square^k}^\bullet|_{\D^{n-i}}\big)$ grow at most polynomially fast.
Therefore, by applying the chain rule to the composition above we get that all partial
derivatives of $\nu$ tend to zero as $v_m\ra 0\in \D^{n-i}$. 
This proves {\bf A($\square^k$)}.\\

Note that the convergence of the derivatives of
$\nu$ to zero shown above is exponentially fast. This together with the fact that
(see item (2) above) the map $\mu_{\square^k}$ is a product on $\D^{n-i}\times\dsquare^k$
imply {\bf B($\square^k$)}. This proves the proposition.\\

\vspace{.6in}

\noindent {\bf \large  Appendix C. Proof of Proposition 4.2.}\\

Recall that $\cA=\Big\{\big( h^\bullet_{\square^k} \,,\, \D^{n-k}\times\dsquare^k
  \big)   \Big\}$ is a normal atlas on $K$, that generates
the normal smooth structure $\cS'$. Also $\{H\0{\square}\}$ is a normal atlas for
$K_X=K_X^{{\mbox{\tiny piece-by-piece}}}$, generating the smooth structure $\cS\0{K\0{X}}$. We will assume that the charts $H\0{\square^k}:
\D^{n-k}\times\dsquare^k\ra K_X$ are defined on the larger sets
$\D^{n-k}(1+\delta)\times\dsquare^k$  (here $\D(1+\delta)$ is the open disc of
radius $1+\delta$). We can obtain this using 2.1.\\

Write $H'_{\square^k}=\Phi\circ H\0{\square^k}$.
It is enough to prove that the maps  $H'_{\square^k}:\D^{n-k}\times\dX\0{\square}
\ra K_X'$ are $C^1$-embeddings.
To prove this we need to prove that the following
coordinate maps are both $C^1$

$$
\begin{array}{llll}
q\0{K}\circ H'\0{\square^k}:&\D^{n-k}\times\dX\0{\square^k}&\longrightarrow&(K,\cS')\\
q\0{X}\circ H'\0{\square^k}:&\D^{n-k}\times\dX\0{\square^k}&\longrightarrow&X.
\end{array}
$$

\noindent  We prove this by induction down the dimension of the skeleta.
First for $k=n$ recall that $H\0{\square^n}:\dX\0{\square^n}\ra K_X$ is just the inclusion. Hence Proposition 3.2 implies that $q\0{K}\circ H'\0{\square^n}:\dX_{\square^n}\ra \dot{\square}^n$ is 
the map $\iota\circ f$, where $\iota:\dot{\square}^n\ra (K,\cS')$ is the inclusion, which is smooth.
(Recall that the inclusion $\square^n\ra (K,\cS')$ is not necessarily differentiable but its restriction $\iota$ to $\dsquare^n$ is smooth; see Remark 1 before Theorem 1.1.1).
Therefore $q\0{K}\circ H'\0{\square^n}$ is smooth.
Also, by the definition of the map $\Phi$, we have  $q\0{X}\circ H'\0{\square^n}=P$, which is also smooth. Moreover, by Proposition 3.1,  $P|\0{\dot{\square}^n}$ is an embedding. Therefore $H'\0{\square^n}$ is a smooth embedding for every $n$-cube $\square^n\in K$.\\

Assume we have proved that  $H'\0{\square^j}$ is a $C^1$-embedding for  every $j$-cube $\square^j\in K$, $j>k$.
We have to prove that the same is true for all $k$-cubes.
We prove this in three parts. In the first part we prove that
$q\0{X}\circ H'\0{\square^k}$ is $C^1$. In the second part we prove that
$q\0{K}\circ H'\0{\square^k}$ is $C^1$. This two parts imply that
$H'\0{\square^k}$ is $C^1$. Finally, in the third part we prove that
$H'\0{\square^k}$ is an embedding.  Fix a $k$-cube $\square^k$.\\

\noindent {\bf FIRST PART.} {\it The map  $q\0{X}\circ H'\0{\square^k}$ is $C^1$.}\\

\noindent {\bf Proof.}
Denote by $V_{\square}$ the image of $H\0{\square}$. For each $\square^i$ with $\square^k<\square^i$ we have that on
$U_{\square^i}=(H\0{\square^k})^{-1}(V_{\square^i})$
we can write $$q\0{X}\circ H_{\square^k}'\,\,=\,\,\Big(q\0{X}\circ\Phi\circ H\0{\square^i}\Big)\circ\Big(H\0{\square^i}^{-1}\circ H\0{\square^k} \Big)
\,\,=\,\,\Big(q\0{X}\circ H'\0{\square^i}\Big)\circ\Big(H\0{\square^i}^{-1}\circ H\0{\square^k} \Big)
$$
\noindent  which is $C^1$ by inductive hypothesis and Proposition
4.1. Since $\Big(\D^{n-k}-\{0\}\Big)\times\dX_{\square^k}$ is contained in the union of the $U_{\square^i}$, $i>k$, we have that
$q\0{X}\circ H'\0{\square^k}$ is $C^1$ outside $\dX_{\square^k}=\{0\}\times\dX_{\square^k}$.\\

Since the map $q\0{X}\circ \Phi|\0{X_{\square^n}}$ can be identified with the map $P:X\ra X$ for an $n$-cube $\square^n$
(recall $X\0{\square^n}$ is a copy of $X$), Proposition 3.1 implies that the derivatives at a point $(0,p)\in \{0\}\times\dX_{\square^k}$
in the $X$ directions $(0,v)$ exist because $q\0{X}\circ H'\0{\square^i}$ 
on $\dX_{\square^k}=\{0\}\times\dX_{\square^k}$ is an embedding.
We next show that the derivatives in the radial directions also exist and vanish.
For this take a ray
$\alpha(t)= (tu,p)\in \D^{n-k}\times\dX_{\square^k}$ and write
$\beta(t)= H\0{\square^k}(\alpha(t))$. 
Note that $\beta'(0)=DH\0{\square^k}.u$ is normal to $X\0{\square^k}$.
We have that
the image of $\beta$ is contained in some $X_{\square^n}$.
As mentioned above the map $q\0{X}\circ \Phi$ on $X_{\square^n}$
can be identified with the map $P:X\ra X$. The fact that the radial derivative
in the direction $u$ exits and vanishes now follows from the addendum to Proposition 3.1.\\

Finally we need to prove that the first derivatives are continuous. But
this follows from a result analogous to Lemma 1.4.1 with
$X_{\square^i}$ replacing $i$-cubes, which can easily be verified.
This concludes the proof of the first part.\\

\noindent {\bf SECOND PART.} {\it The map  $q\0{K}\circ H'\0{\square^k}:
\D^{n-k}\times\dX_{\square^k}\ra (K,\cS')$ is $C^1$.}\\

\noindent {\bf Proof.} This proof will take the next five pages.
First note that, by Corollary 3.3 and the definition
of $H\0{\square}$ we have

\begin{equation*} q\0{K}\circ H'\0{\square^i}(\,t\,v\,,\, p\,)=F\,\Big( exp\0{p}\,\big(\,2\,r\,t\,\,h_{\square^i}(v)\,\big)\,\Big).
\tag{1}
\end{equation*}

\noindent Write $G_{\square^i}=\big(h^\bullet_{\square^i}\big)^{-1}\circ q\0{K}\circ H'_{\square^i}:\D^{n-i}\times\dX_{\square^i}\longrightarrow\D^{n-i}\times\dsquare^i$. Since $\{h^\bullet\0{\square}\}$ is an atlas for $(K,\cS')$,
by inductive hypothesis we have that $G_{\square^i}$ is $C^1$,
for $i>k$, and \\

\noindent {\bf (C.1.)}  \,{\it to prove that $q\0{K}\circ H'\0{\square^k}$ is $C^1$
 it is enough to prove that $G_{\square^k}$ is \, $C^1$.}\\

\noindent Write also $G_{\square}=(R_\square\, ,\, T_\square)$. For
$u=tv\in \D^{n-i}$, $t=|u|$, we have

\begin{equation*}
G_{\square^i}(\,t\,v\,,\, p\,)\,\,=\,\,\Big(R_{\square^i}(\,t\,v\,,\, p\,)\, ,\, T_{\square^i}(\,t\,v\,,\, p\,)  \Big)\,\,=\,\, \big(h^\bullet_{\square^i}\big)^{-1}
\circ F\,\Big( exp\0{p}\,\big(\,2\,r\,t\,\,h_{\square^i}(v)\,\big)\,\Big)
\tag{2}
\end{equation*}

\noindent 
It follows from (2) and Lemma 2.6 (iii), (iv), that we can write

\begin{equation*}
G_{\square^i}(\,u\,,\, p\,)\,\,=\,\,\Big(R_{\square^i}(\,u\,)\, ,\, T_{\square^i}(\,|u|\,,\, p\,)  \Big)\tag{3}
\end{equation*}

\noindent  that is, $R$ does not depend on $p$ and $T$ depends on $p$ and the
length $|u|$ of $u$ (not on the direction of $u$). Also it can be checked from
Lemma 2.6 (iv) and Lemma 2.7 that $(u, p)\mapsto T_{\square}(|u|, p)$ is smooth.
This together with (C.1) and (3) imply that\\

\noindent {\bf (C.2.)} {\it to prove that $q\0{K}\circ H'\0{\square^k}$ is $C^1$ it is enough to prove that
$R_{\square^k}:\D^{n-k}\ra\D^{n-k}$ is \, $C^1$.}\\

\noindent {\bf Claim C.3.} {\it The map\, $R_{\square^k}$ is $C^1$ on $\D^{n-k}-\{ 0\}$.}\\

\noindent {\bf Proof of Claim C.3.} Recall that by inductive hypothesis we have that $G_{\square^i}$, $R_{\square^i}$ and $q\0{K}\circ H\0{\square^i}'$ are $C^1$, for all $\square^i$, $i>k$.
Denote by $V_{\square}$ the image of $H\0{\square}$. For each $\square^i$ with $\square^k<\square^i$ we have that on
$U_{\square^i}=(H\0{\square^k})^{-1}(V_{\square^i})$
we can write $$q\0{K}\circ H_{\square^k}'\,\,=\,\,\Big(q\0{k}\circ\Phi\circ H\0{\square^i}\Big)\circ\Big(H\0{\square^i}^{-1}\circ H\0{\square^k} \Big)
\,\,=\,\,\Big(q\0{k}\circ H'\0{\square^i}\Big)\circ\Big(H\0{\square^i}^{-1}\circ H\0{\square^k} \Big)
$$
\noindent  which is $C^1$ by inductive hypothesis and Proposition
4.1. Since $\D^{n-k}-\{0\}=(\D^{n-k}-\{0\})\times\{p\}$ (for any $p\in
\dX\0{\square^k}$) is contained in the union of the $U_{\square^i}$, $i>k$, we have that
$q\0{K}\circ H'\0{\square^k}$ (hence $G_{\square^k}$ and $R_{\square^k}$) is $C^1$ outside $0$.  This proves Claim C.3.\\

From  Lemma 2.5 (ii) and the fact that the derivative of the exponential (at 0)
is the identity we get

\begin{equation*}
\frac{\p}{\p v}R_{\square^k}(0)=v.
\tag{4}
\end{equation*}

\noindent That is, all directional derivatives at 0 of $R_{\square^k}$ exist and
if $R_{\square^k}$ were differentiable its derivative at 0 would be the identity matrix 1.
It follows from (4), (C.2) and Claim C.3 that\\

\noindent {\bf (C.4.)} {\it to prove that $q\0{K}\circ H'\0{\square^k}$ is $C^1$ it suffices to prove
$DR_{\square^k}|\0{q}\ra 1$ (the identity matrix)  

\hspace{.21in} as $q\ra 0$.}\\

Write $S=\sL(\square^k,K)=\sL(X_{\square^k}, K_X)$ (at some point $F(p)\in\dsquare^k$ and
$p\in\dX\0{\square^k}$, respectively, and recall we are using ``direction" links).
Also write $\D^{n-k}=\D^{n-k}\times \{p\}\sbs \D^{n-k}\times \dX_{\square^k}$.
For $\square^j$, $\square^k\sbs\square^j$ set $\sigma\0{\square^j}=
\D^{n-k}\cap(H_{\square^k})^{-1}(X\0{\square^j})$ and
$\dsigma\0{\square^j}=
\D^{n-k}\cap(H_{\square^K})^{-1}(\dX\0{\square^j})$.
Note that the sets $\sigma\0{\square^j}$ and $\dsigma\0{\square^j}$ are cone sets.
That is, if $u\in\sigma\0{\square^j}$ then $tu\in\sigma\0{\square^j}$, $t\in [0,1]$.
Similarly for $\dsigma\0{\square^j}$ (see Lemma 2.6).\\

\noindent {\bf Claim C.5.} {\it We have 
 $R_{\square^k}(\sigma\0{\square^j})=\sigma\0{\square^j}$
\,\,and\,\, $R_{\square^k}(\dsigma\0{\square^j})=\dsigma\0{\square^j}$.}\\

\noindent {\bf Proof of Claim C.5.} 
We prove the first identity, the second one is similar.
Let $u=tv\in \dsigma\0{\square^j}$, $|v|=1$. We assume $t>0$.
Then $exp\0{p}(2\,r\,t\,h\0{\square^k}(v))=H\0{\square^k}(u,p)\in X\0{\square^j}$. 
By (1) and Lemma 2.6 (ii) we have that $h^\bullet\0{\square^k}\circ \,G\0{\square^k}(u,p)\in \square^j$. By the definition of  $h^\bullet\0{\square^k}$ we get
$\big(a\,h\0{\square^k}(\frac{1}{a}\,R_{\square^k}(u))\,,\,T_{\square^k}(|u|,p)\,\big)\in \square^j$, where $a$ is the length of $R_{\square^k}(u)$.
By Lemma 2.6 (ii) we get  $h\0{\square^k}(\frac{1}{a}\,R_{\square^k}(u))\in T_pX\0{\square^j}$, which implies $\frac{1}{a}\,R_{\square^k}(u)\in\sigma\0{\square^j}$. Since this set is a cone set it follows that
$R_{\square^k}(u)\in\sigma\0{\square^j}$. This proves Claim C.5.\\

Now, let $(q_m)$ be a sequence in $\D^{n-k}$, with $q_m\ra0$, as $m\ra \infty$. We can assume (arguing by contradiction) that 
$q_m= t_m u_m$, with $(t_m,u_m)\in\R^+ \times\bS^{n-k-1}$, $t_m\ra 0$, $u_m\ra u\in \dsigma\0{\square^j}$
for some $\square^j\in K$ containing $\square^k$. Hence:\\

\noindent{\bf  (C.6.)} {\it to prove that $q\0{K}\circ H'\0{\square^k}$ is $C^1$ it suffices to prove that
$DR_{\square^k}|\0{(t_mu_m)}\ra 1$, as $m\ra\infty$, 

\hspace{.26in}where
$u_m\ra u\in\dsigma\0{\square^j}$ and $t_m\ra 0$.}\\

\noindent {\bf Claim C.7.} {\it Statement (C.6) holds for  $j=n$.}\\
 
\noindent {\bf Proof of Claim C.7.} We have that $u\in\dsigma\0{\square^n}$ for some $\square^n\in K$. Therefore there is a small compact neighborhood $V$
of $u\in\bS^{n-k-1}\cap \dsigma\0{\square^n}\sbs\D^{n-k}$ such that (we can assume that) all $q_m$ and $tu$, $t\in (0,1]$, 
lie in the interior of the cone $\rC V$. Denote by $h:\D^{n-k}\ra\rC S$ the map $\rC h_{\square^k}$, where
$h_{\square^k}$ is the link smoothing of $S=\sL(\square^k,K)$. Since $F|\0{X_{\square^n}}=f$, on $\rC V$ we can write
$R_{\square^k}=h^{-1}\circ (\pi\circ f\circ e)\circ h$, where $e$ is the exponential map given by
$e(v)=exp\0{p}(2r v)$, and $\pi:\rC S\times \square^k\ra \rC S$ is the
projection. (Note that $e(tv)=E(2rt,v)$, where $E$ is as in Lemma 2.6.) Hence 
$(DR_{\square^k})|\0{q_m}=(Dh|\0{y_m})^{-1}\, D(\pi\circ f\circ e)|\0{h(q_m)}\, Dh|\0{q_m}$, where 
$y_m=h^{-1}(\pi\circ f\circ e)(h(q_m))$. By Lemma 2.5 (ii) and the fact that
the derivative of the exponential at 0 is the identity we have that $ D(\pi\circ f\circ e)|\0{h(q_m)}\ra 1$, as $q_m\ra 0$.
On the other hand, since $h$ is a cone map, by Lemma 1.4.2 we get that  $Dh$ and $Dh^{-1}$ are both bounded on $\rC V$. Moreover,
$Dh|\0{q_m}=Dh|\0{u_m}$ (see remark after Lemma 1.4.2) and $Dh|\0{y_m}=Dh|\0{\frac{y_m}{|y_m|}}$.
But, since $\pi\circ f\circ e$ is smooth and $D(\pi\circ f\circ e)|\0{p}=1$ we have that for any $v$ we get

{\small \begin{equation*}
\frac{\pi\circ f\circ e\,\big(\, t_m\, v\big)}{t_m}\,\,\longrightarrow\,\, v.
\tag{5}
\end{equation*}}

\noindent From the fact that $h^{-1}$ is a cone map and (5) we have

{\small $$\frac{t_m}{|\,\,\,h^{-1}\Big((\pi\circ f\circ e)\,\big(t_m\,h\, (  u_m) \big)  \Big)\,\,\,  |}=
\Big[\,\,|\,h^{-1}\Big(\frac{(\pi\circ f\circ e)\,\big(t_m\,h\, (  u_m)\big) }{t_m} \,\, \Big)\,\,|\,\,\Big]^{-1}\,\,\longrightarrow\,
lim\0{m\ra\infty}\,|u_m|^{-1}\,\,=\,\,1 $$}

\noindent  This together with (5) and the fact that $h$ and $h^{-1}$ are cone maps imply
{\small $$\begin{array}{lll}
\frac{y_m}{|y_m|}&=&\frac{h^{-1}\Big((\pi\circ f\circ e)\,h\, \big( t_m\, u_m \big)  \Big)}
{|\,\,\,h^{-1}\Big((\pi\circ f\circ e)\,h\, \big( t_m\, u_m \big)  \Big) \,\,\, |}\,\,=\,\,
h^{-1}\bigg(\,\,
\frac{\Big((\pi\circ f\circ e)\,h\, \big( t_m\, u_m \big)  \Big)}
{|\,\,\,h^{-1}\Big((\pi\circ f\circ e)\,h\, \big( t_m\, u_m \big)  \Big)\,\,\,  |} \,\, \bigg)\\\\
&=&h^{-1}\bigg(\,\,
\frac{\Big((\pi\circ f\circ e)\,h\, \big( t_m\, u_m \big)  \Big)}
{t_m} \frac{t_m}{|\,\,\,h^{-1}\Big((\pi\circ f\circ e)\,\big(t_m\,h\, (  u_m) \big)  \Big)\,\,\,  |} \,\,\bigg)\,\,\xrightarrow{\,\,\,\,\,\,\,\,}\,\,u.
\end{array}$$}

\noindent  Consequently
$\frac{y_m}{|y_m|}\ra u$. Therefore $DR_{\square^k}|\0{q_m}\ra 1$. This proves Claim C.7.\\

We will prove statement (C.6) by decreasing induction on $j$. Claim C.7 was
the first step of this induction. We assume statement C.6 holds for all
$\square^i$, $j<i$.\\

\noindent {\bf Claim C.8.} {\it Statement (C.6) holds for  $j$.}\\
  
\noindent {\bf Proof of Claim C.8.} 
The proof has two steps.\\

\noindent {\it Step 1. It is enough to assume that $u_m=u$. }\\
As in the proof of Claim C.7 let $V$ be 
a small compact neighborhood
of $u\in\bS^{n-k-1}\sbs\D^{n-k}$ such that all $q_m$ and $tu$, $t\in (0,1]$, 
lie in the interior of the cone $\rC V$. From the definition of $G_\square$
we have that on $\rC V$ we can  write

$$
G_{\square^k}\,\,=\,\, \Big( \,\big( h^\bullet_{\square^j}\big)^{-1}\circ\,\, h^\bullet_{\square^k}   \Big)^{-1}\,\,
\circ\,\, G_{\square^j}\,\,\circ\,\, \Big( \,\big( H_{\square^j}\big)^{-1}\circ H_{\square^k}   \Big).
$$

To simplify the notation write 
$h= \big( h^\bullet_{\square^j}\big)^{-1}\circ\, h^\bullet_{\square^k}$ and $H=\big( H_{\square^j}\big)^{-1}\circ H_{\square^k}$. Hence 

\begin{equation*}DG_{\square^k}=Dh^{-1}.DG_{\square^j}.DH.
\tag{6}
\end{equation*}

We next compare $DG_{\square^k}|\0{(t\0{m} u\0{m},p)}$ and $DG_{\square^k}|\0{(t\0{m} u,p)}$. 
We analyze the three terms $DH$, $DG_{\square^j}$, $Dh$ in (6).\\

\noindent {\bf First term:  $DH$.} Since $H$ is a cone map we get that
$DH|\0{t\0{m}u_{m}}-DH|\0{t\0{m}u}\ra 0$.\\

\noindent {\bf Remark.} The map $H$ is a Euclidean-to-hyperbolic cone map, and it
is not a euclidean cone map but it is a euclidean
cone map up to a smooth change of coordinates on a compact set.\\

\noindent {\bf Second term:  $DG_{\square^j}.$}
Differentiating (3) we get

\begin{equation*} DG_{\square^j}|\0{(u,y)} (v,w)=\Big(\,  DR_{\square^j}|\0{u}\,.\,v\,\,,\,\, 
\frac{\p}{\p t}T_{\square^j}|\0{(t,y)}\frac{u.v}{|u|}\,+\, \frac{\p}{\p y}T_{\square^j}|\0{(t,y)} .w  \,\Big)
\tag{7}
\end{equation*}

\noindent where $t=|u|$. It can be
checked from Lemma 2.6(iv) (see also (3)) that $T_{\square^j}$ can be extended to a smooth map
on $\bar{\D}^{n-k}\times X_{\square^j}$ (which is compact) 
the second term (i.e the $T_\square$ term) is  Lipschitz on the variables
$t$ and $y$. Since the distance between 
$H(t_mu_m,p)$ and $H(t_m u,p)$   goes to zero it follows that the $T_\square$ terms
in the right hand side of the equation above evaluated at $H(t_mu_m,p)$ and $H(t_m u,p)$
get close as $m\ra\infty$. Also, by inductive hypothesis, the first terms tend both to 1.
Therefore we get
$$
DG_{\square}|\0{H(t\0{m}u\0{m})}-DG_{\square}|\0{H(t\0{m}u)}\longrightarrow\,\,
\,\,0
$$

\noindent as $m\ra \infty$. \\

\noindent {\bf Third term: $Dh$.}
Since $h$ is a cone map to prove that $Dh\0{G_{\square^j}(H(t\0{m}u\0{m}))}$
and $Dh\0{G\0{\square^j}(H(t\0{m}u))}$ are close we need to prove that
the directions of
$G_{\square^j}(H(t\0{m}u\0{m}))$ and $G_{\square^j}(H(t\0{m}u))$
 are close. This is equivalent to proving that the directions of their images by $h$ are close. Since 
 $G_{\square^k}=h\circ G_{\square^j}\circ H$ this means
proving that the directions of $G_{\square^k}(t\0{m}u\0{m})$ and $G_{\square^k}(t\0{m}u)$ are close.  
Let $E$ and $p\0{l}$ be the maps in Lemma 2.6.
We can assume (arguing by contradiction) that $t_mu_m$, $t_m u$ lie on
$ X_{\square^n}$, for some $\square^n$.
Since $h^\bullet_{\square^k}$ is also
a cone map (on the first variable) it is enough to prove that the directions of $p\0{l}\circ f\circ E(2r t_m, h_{\square^k}(u_m))$
and $p\0{l}\circ f\circ E(2rt_m, h_{\square^k}(u))$ are close.
But this is true because $p\0{l}\circ f\circ E$ is smooth.
This concludes step 1.\\

\noindent{\it Step 2. We prove that
 $DR_{\square^k}|\0{(t_mu)}\ra 1$, as $m\ra\infty$, where
$ u\in\dsigma\0{\square^j}$ and $t_m\ra 0$.}\\

Note that
every inclusion $\dsigma\0{\square^j}\hookrightarrow \D^{n-k}$ is a smooth embedding (see Section 1). We have two cases.\\

\noindent {\bf First case.} {\it We have that
 $DR_{\square^k}|\0{(t_mu)}v\ra v$, as $m\ra\infty$, when $v$ is tangent to $\sigma\0{\square^j}$.}\\

\noindent This follows from an argument similar to the one given
in the proof of (C.7) (recall that from (C.5) we have  $R_{\square^k}(\sigma\0{\square^j})=\sigma\0{\square^j}$). This proves the first case.\\\\

Recall that $H$ is the change of variables
$H=\big( H_{\square^j}\big)^{-1}\circ H_{\square^k}$.
Let $u\in\dsigma\0{\square^j}$ and $v\in\R^{n-k}=T_u \D^{n-k}$.
% we have that $H(tu)\in \{ 0\}\times\dX_{\square^j}$ and write $H(tu)=(0, q\0{t})$. %The change of charts $H^{-1}$ gives a fiber bundle over $c\sigma\0{\square^j}$ in %$\D^{n-k}$: the fibers are the images
%of $\D^{n-j}\times\{x\}\sbs\D^{n-j}\times \dX_{\square^j}$ by $H^{-1}$, where $x$ %is in the image of $c\sigma\0{\square^j}$ by $H$.
%Note that these fibers are transversal to $\dsigma\0{\square^j}$.
%If $v$ is tangent to  one of these fibers, at some point $u\in\dsigma\0{\square^j}$
We say that $v$ is an $X$-{\it fiber vector at $u$} if
$DH|\0{u}v=(z,0)\in
\R^{n-j}\times T_{H(u)}\dX_{\square^j}=T_{H(u)}(\D^{n-j}\times\dX_{\square^j})$, for some $z\in\R^{n-j}$. We write $v=v\0{z}^X$
(the reason for the upper index $X$ will be clear in a moment). Fixing $z$ we
obtain a constant vector field $(z,0)$, hence we obtain the corresponding
vector field $v\0{z}^X$ of $X$-fiber vectors on $\dsigma\0{\square^j}$. Thus 
$v\0{z}^X$ is characterized by $DH|\0{u}.v\0{z}^X(u)=(z,0)$. Equivalently $v\0{z}^X(u)=\big(DH|\0{u}\big)^{-1}.(z,0).$\\

Similarly, we can work on $K$ instead of $K_X$, and $h$ instead of $H$
and obtain vector fields of $\square$-{\it fiber vectors} $v\0{z}^\square$
on $\dsigma\0{\square^j}$ with the property that $Dh|\0{u}.v\0{z}^\square(u)=(z,0)
\in\R^{n-j}\times T_{h(u)}\square^j=T_{h(u)}(\D^{n-j}\times\square^j)$.  Equivalently 
\begin{equation}
v\0{z}^\square(u)=\big(Dh|\0{u}\big)^{-1}.(z,0).
\tag{8}
\end{equation}

\noindent {\bf Claim C.9.} {\it We have $v\0{z}^X=v\0{z}^\square$.}\\

\noindent {\bf Proof of Claim C9.} Fix $u$. Then the (hyperbolic) geodesic
$t\mapsto H\0{\square^k}(tu)$, and the straight segment 
$t\mapsto h\0{\square^k}^\bullet(tu)$ are both contained in $X\0{\square^n}$
and $\square^n$, respectively, for some $\square^n$. This together with
the fact that both $h\0{\square^k}^\bullet$ and $H\0{\square^k}$ use
the same link smoothing $h\0{\square^k}$ in their definition imply that
we can reduce our problem to the following setting. Consider $\R^n=\R^{n-k}\times
\R^k$, $\R^k\sbs\R^j$, with canonical
metrics $\sigma\0{\R^n}$ and $\sigma\0{\HH^n}=
\sigma\0{\HH^{n-k}}+cosh^2(r)\sigma\0{\HH^{k}}$, and $z$ (a constant vector field) perpendicular to $\R^j$. Here $r$ is the distance to $\HH^k$. In this case $h\0{\square^k}^\bullet$
corresponds to the perpendicular (to $\R^k$) exponential map from a point $p\in \R^k$ 
and $H\0{\square^k}$ corresponds to
the perpendicular (to $\HH^k=(\R^k,\sigma\0{\HH^k})$) exponential map from $p$. 
The former exponential is just the inclusion and the latter exponential is done with respect to $\sigma\0{\HH^n}$. But in this setting
these two exponentials coincide, hence the preimage of $z$ by them also
coincide. This proves the claim.\\

Given $z$ as above we write $v\0{z}$  to denote $v\0{z}^X$ and $v\0{z}^\square$
and we say that $v=v\0{z}$ is  a {\it fiber vector}.
The following statement can be easily verified in the Euclidean case (i.e for
$v\0{z}^\square$).

\begin{equation*}
v\0{z}(tu)\,\,=\,\,v\0{z}(u).
\tag{9}
\end{equation*}\vspace{.2in}

To prove Step 2 it is enough to prove the following.\\

\noindent {\bf Second case.} {\it 
If $v$ is fiber vector at  $u$, then 
$DR_{\square^k}|\0{(t\0{m}u)}v\ra v$, as $m\ra\infty$.}\\

\noindent We have $v=v\0{z}(u)$, hence 
\begin{equation}DH|\0{(t\0{m} u)}.v=(z,0).
\tag{10}
\end{equation}

Using formula
(7), the inductive hypothesis and the fact that $(0,q\0{t})=H(tu)\in \{ 0\}\times\dX_{\square^j}$
we see that  for any $z'\in\R^{n-j}$ we have
$$
DG_{\square^j}|\0{(H(t_m u))}(z',0)=(z',0).
$$ 
\noindent This together with the Equation (6) imply
\begin{equation}
DG_{\square^k}|\0{(t_mu)}(v,0)=Dh^{-1}\0{G\0{\square^j}(H(t_mu))}.DH\0{t\0{m} u}.v.
\tag{11}
\end{equation}

\noindent Therefore, substituting  (8), (10) and Claim 9
into Equation (11), we get
\begin{equation}
DG_{\square^k}|\0{(t_mu)}(v,0)=DG_{\square^k}|\0{(t_mu)}(v\0{z}(u),0)=
v\0{z}\big( G_{\square^k}(t_mu)\big).
\tag{12}
\end{equation}

%THIS IS OK, IT IS $k$ AND NOT $j$ IN THE LAS INDEX ABOVE!!!

\noindent  Note that, by (7) and Remark 2.7,
$DG_{\square^k}|\0{(t_mu)}(v,0)=\big( DR_{\square^k}|\0{(t_mu)}.v,0)$. This together with (12) imply that to prove Case 2 we have to prove that $v\0{z}\big( G_{\square^k}(t_mu)\big)=(v,0)=v_z(u)$. 
Consequently,
since $h$ is a cone map it is enough to prove that the directions of $G_{\square^k}(t_mu)=h^{-1}\circ G_{\square^j}\circ H(t_mu)$ tend 
to $u$, as $m\ra\infty$. 
But this last statement is implied by $$lim_{m\ra\infty}\frac{G_{\square^k}(t_mu)}{t_m}=DG_{\square^k}|\0{(0)}.u=\Big( DR_{\square^k}|\0{0}.u\,,\, 0 \Big)=u$$

\noindent which follows from equation (4). This proves the second case, Step 2, Claim C.8 and concludes the second part.\\

\vspace{.2in}

\noindent {\bf THIRD PART.} {\it The maps $H'_{\square}$ are $C^1$-embeddings.}\\

\noindent{\bf Proof.} Again by induction. This is true for
$H'_{\square^n}:\dX_{\square^n}\hookrightarrow K'_X$, which can be
identified with $P|\0{\dsquare^n}$ (see Proposition 3.1).
Assume that  the $H'_{\square^j}$ are $C^1$-embeddings for $j>k$. Fix a $\square^k$.
Using the argument used in the first and second parts we get that 
$H'_{\square^k}$ is a $C^1$-embedding outside $\dX_{\square^k}
=\{0\}\times\dX_{\square^k}$.
From the second part we see that the derivative $DH'_{\square^k}$
maps non-zero vectors $v$ at $\dX_{\square^k}$ in the $\D^{n-k}$ direction to non-zero vectors (see Equation (4)).
As mentioned in the first part the map $q\0{X}\circ H'_{\square^k}|\0{\dX_{\square^k}}$ can be identified with $P|\0{\dX\0{\square^k}}$, which
is a diffeomorphism (see Proposition 3.1). Hence $DH'_{\square^k}$ maps
non-zero vectors $w$ in the $X_{\square^k}$ direction to non-zero vectors.
Moreover, $DH'_{\square^k}.v$ is perpendicular to $DH'_{\square^k}.w$.
It follows that $H'_{\square^k}$ is an embedding on $\{0\}\times\dX_{\square^k}$.
This concludes the third part and completes the proof of Proposition 4.2.

Pedro Ontaneda

SUNY, Binghamton, N.Y., 13902, U.S.A.

\end{document}